\documentclass[12pt]{amsart}

\usepackage{amsmath}
\usepackage{amssymb}
\usepackage{fancybox}
\usepackage{eucal}
\usepackage{amsthm}
\usepackage{amscd}
\usepackage{wasysym}
\usepackage{url}
\usepackage{MnSymbol} %symbols for infinitary formulas

\usepackage[all,cmtip]{xy}

\usepackage{enumitem}% http://ctan.org/pkg/enumitem

\usepackage{hyperref}

\newtheorem{thm}{Theorem}[section]
\newtheorem{prop}[thm]{Proposition}

\newtheorem{lemma}[thm]{Lemma}

\newtheorem{cor}[thm]{Corollary}

\newtheorem{definitiontemp}[thm]{Definition}
\newenvironment{defn}{\begin{definitiontemp} 
\normalfont}{\end{definitiontemp}}

\def\oCK{\omega_1^{\rm CK}}

\newenvironment{pf}{\begin{trivlist}\item[\hskip\labelsep
{\it Proof.}]}{\end{trivlist}}

\newcommand{\comment}[1]{}

\newcommand{\set}[2]{\ensuremath{ \{ #1 : #2 \} }}

\newcommand{\N}{\mathbb{N}}
\newcommand{\Z}{\mathbb{Z}}

\newcommand{\Q}{\mathbb{Q}}
\newcommand{\R}{\mathbb{R}}

\newcommand{\A}{\mathcal{A}}
\newcommand{\B}{\mathcal{B}}
\newcommand{\D}{\mathcal{D}}
\newcommand{\U}{\mathcal{U}}

\renewcommand{\P}{\mathcal{P}}

\newcommand{\ep}{\varepsilon}

\newcommand{\la}{\langle}
\newcommand{\ra}{\rangle}

\def\diverges{\!\uparrow}
\def\converges{\!\downarrow}
\newcommand{\at}{\char'100}

\renewcommand{\a}{\text{$\bf{a}$}}
\renewcommand{\b}{\text{$\bf{b}$}}
\renewcommand{\c}{\text{$\bf{c}$}}
\renewcommand{\r}{\text{$\bf{r}$}}

\newcommand{\f}{\text{$\bf{f}$}}
\newcommand{\g}{\text{$\bf{g}$}}
\newcommand{\h}{\text{$\bf{h}$}}
\newcommand{\x}{\text{$\bf{x}$}}
\newcommand{\y}{\text{$\bf{y}$}}

\newcommand{\dom}[1]{\text{dom}(#1)}

\def\s01{\ensuremath{\Sigma^0_1}}
\def\d02{\ensuremath{\Delta^0_2}}
\def\phi{\varphi}

\def\res{\!\!\upharpoonright\!}

\begin{document}

\title{Effectivizing Lusin's Theorem}

\author[R. Miller]{Russell Miller}
\address{Dept.\ of Mathematics, Queens College, \& Ph.D.\ Programs in Mathematics \& Computer Science,
Graduate Center, City University of New York,  USA}
\email{Russell.Miller@qc.cuny.edu}
\urladdr{\href{http://qcpages.qc.cuny.edu/~rmiller/}{http://qcpages.qc.cuny.edu/$\sim$rmiller/}}

\thanks{
The author was partially supported by Grant \#DMS -- 1362206 from
the National Science Foundation, Grant \#581896
from the Simons Foundation, and several grants from
the PSC-CUNY Research Award Program.  This work was initiated
at a workshop held at the Schlo\ss\ Dagstuhl
in February 2017, where Vasco Brattka aided its development
with useful conversations.  It was subsequently improved
by several suggestions from an anonymous referee.}

\begin{abstract}
Lusin's Theorem states that, for every Borel-measurable function $\f$
on $\R$ and every $\ep>0$, there exists a continuous function $\g$ on $\R$
which is equal to $\f$ except on a set of measure $<\ep$.  We give a proof
of this result using computability theory, relating it to the near-uniformity
of the Turing jump operator, and use this proof to derive several uniform
computable versions.  Easier results, which we prove by the same methods,
include versions of Lusin's Theorem with Baire category in place of
Lebesgue measure and also with Cantor space $2^\N$ in place of $\R$.
The distinct processes showing generalized lowness for generic sets and
for a set of full measure are seen to explain the differences between
versions of Lusin's Theorem.
\end{abstract}

\maketitle

\section{Introduction}
\label{sec:intro}

Lusin's Theorem, first proven by Nikolay Lusin (or Luzin) in 1912,
states informally that every measurable function on the real numbers $\R$
is ``nearly continuous,'' in terms of the Lebesgue measure $\mu$ on $\R$.
The standard formal statement is as follows.

\begin{thm}[Lusin's Theorem, 1912]
%\label{thm:Lusin}
For every Borel-measurable function $\f:\R\to\R$ and every $\ep>0$,
there exists a continuous function $\g:\R\to\R$ such that
$$ \mu(\set{\x\in\R}{\f(\x)\neq \g(\x)})<\ep.$$
\end{thm}
Alternative versions allow $\pm\infty$ as values of the functions in question.
(Some versions only state that $\f$ restricts to a continuous function
on a set of comeasure $<\ep$, but we will consider the stronger version.)
A common method of proof of this result involves Egorov's Theorem,
that continuity of a function on a compact set is ``nearly'' uniform continuity.
Indeed, in the standard text \cite{Royden}, Lusin's Theorem is posed as an
exercise, following the exposition of Egorov's Theorem.  On the other hand,
\cite[Thm.\ 2.24]{Rudin} proves Lusin's Theorem from basic principles
(and later poses Egorov's Theorem as an exercise).

The purpose of this article is to show that Lusin's Theorem is a direct consequence
of the fact, well known in computability theory, that for almost all sets $X$
of natural numbers, the jump $X'$ of $X$ satisfies the property
$$ X' \leq_T \emptyset'\oplus X,$$
i.e., the jump, which is essentially the Halting Problem relativized to $X$,
can be computed merely from the Halting Problem $\emptyset'$
along with $X$ itself.  (The \emph{join}
$\emptyset'\oplus X=\set{2n}{n\in\emptyset'}\cup\set{2n+1}{n\in X}$
is a least upper bound of $\emptyset'$ and $X$ under $\leq_T$.)
Certain sets $X$ fail to satisfy this property,
but those sets form a meager subclass of measure $0$ within the Cantor space $2^\N$.
The computability-theoretic reason why Lusin's Theorem needs to use an $\ep>0$
is that the Turing reduction from $X'$ to $\emptyset'\oplus X$ is only uniform
up to a set of measure $\ep$.  In contrast, when we use Baire category theory,
comeager-many $X\in 2^\N$ satisfy $X'\leq_T\emptyset'\oplus X$,
and here a single Turing reduction succeeds uniformly on the entire class.
This accounts for the companion theorem in analysis, which is much easier to prove
and requires no $\ep$-fudging:
that every Borel function $\f$ restricts to a continuous function whose domain is
a comeager subset of $\R$.

Our proof of Lusin's Theorem demands more background
in computability and descriptive set theory than the standard proofs.
The point is not to replace the original proofs, but rather to highlight
the connections between Lusin's Theorem and computability theory.
When given in full, our proof will require substantial attention
to detail, but it can be summarized very neatly as the following sequence of steps.
\begin{enumerate}
\item
As shown in \cite{W00}, the continuous functions $\g:\R\to\R$ are precisely those
for which there is a Turing functional $\Gamma$ and an oracle $S\subseteq\N$
such that, for all Cauchy sequences $X$ converging fast to any $\x$,
$\Gamma^{S\oplus X}$ converges fast to $\g(\x)$.
\item
Borel-measurable functions $\f(\x)$ are those which can be described
(in computable analysis) by the action of a Turing functional $\Phi$ whose
oracle is the $\alpha$-th jump of the input $\x$, for some countable ordinal
$\alpha$ that depends on $\f$, along with an oracle set $S\subseteq\N$.
That is, if $X$ is a representation (normally by a fast-converging Cauchy
sequence) of $\x$, then
$$ \Phi^{\left( (S\oplus X)^{(\alpha)}\right)}$$
computes a Cauchy sequence which converges fast to $\f(\x)$.
The ordinal $\alpha$ gives the level of $\f$ in the Borel hierarchy
(sometimes denoted by $F_\sigma$, $G_{\delta\sigma}$, etc.;
also by $\Sigma_\alpha$ or $\Pi_\alpha$).
This generalizes the preceding item, where $\alpha=0$,
\item
For each fixed ordinal $\alpha<\omega_1$, almost all subsets $X\subseteq\N$ have the property
that the $\alpha$-th jump $X^{(\alpha)}$ is Turing-reducible to $\emptyset^{(\alpha)}\oplus X$.
(For $\alpha=1$, such an $X$ is said to be \emph{generalized low$_1$}, or $GL_1$;
this situation was described above.)
More generally, relativizing to a fixed $S\subseteq\N$, almost all $X$
satisfy $(S\oplus X)^{(\alpha)}\leq_T S^{(\alpha)}\oplus X$.  Details appear in \cite{S72}.
\item
In the preceding item, ``almost all'' refers to Lebesgue measure
on the space of inputs $X$, and the Turing reducibility, while not uniform,
is uniform up to a set of arbitrarily small measure, and moreover
is uniform in the bound on that set's measure.  That is,
there is a Turing functional $\Psi$ such that, for every $\ep>0$ and every $S$,
$$ \mu(\set{X}{\Psi^{S^{(\alpha)}\oplus X}(\ep,~\cdot~) \neq (S\oplus X)^{(\alpha)}})<\ep.$$
We often write $\Psi_{\ep}$ for the unary functional $\Psi(\ep,~\cdot~)$.
$\Psi_{\ep}$ is not uniform in $\alpha$.
\end{enumerate}
Taken together, these four items suggest a natural approach for proving Lusin's Theorem:
if $\f(\x)$ is given by $\Phi^{(S\oplus X)^{(\alpha)}}$ as in (2), and $(S\oplus X)^{(\alpha)}$ is given in turn
by $\Psi_{\ep}^{S^{(\alpha)}\oplus X}$, then composing these two should yield
a function
$$ \Gamma^{S^{(\alpha)}\oplus X}= \Phi^{\left(\Psi_{\ep}^{S^{(\alpha)}\oplus X}\right)}$$
which will be continuous by virtue of item (1) and will equal $\f$ up to a set of measure $<\ep$
by virtue of item (4).

Of course, many readers will have spotted potential flaws in this argument already,
and it will require a good deal of work to address them.  In particular, even for
those $X$ on which $\Psi^{\emptyset^{(\alpha)}\oplus X}$ fails to output $X^{(\alpha)}$,
$\Gamma$ must still give a coherent output, i.e., a fast-converging Cauchy sequence.
Moreover, the sequences computed by $\Gamma$ for two distinct $X_0$ and $X_1$
must converge to the same real number whenever $X_0$ and $X_1$ both converge
to the same input $\x$.  Of course, each $\x$ is the limit of continuum-many
fast-converging Cauchy sequences, many of which may not be handled correctly
by $\Psi$, so this appears to be a serious problem.  Finally, Lebesgue measure
on Cantor space $2^{\N}$ is used for the statement that almost all $X$ are
$GL_1$, whereas we must use Lebesgue measure on $\R$, represented as a quotient
of the space of all fast-converging Cauchy sequences, to address Lusin's Theorem.

Nevertheless, we will
overcome these difficulties and give a proof of Lusin's Theorem essentially
following the outline given above.  This will require a reasonable background in
computable analysis and also in the computability-theoretic notion of the jump
operation on subsets of $\N$ and its iteration through the countable ordinals.
Two appendices (Sections \ref{sec:computableanalysis} and
\ref{sec:ordinaljumps}) are devoted to these two topics, as few readers can be
expected to be closely familiar with both of them.  For more extensive presentations,
we suggest \cite{W00} for the computable analysis; \cite{Kechris} for descriptive set theory
such as the preceding characterization of the Borel functions; \cite[Chapter III]{S87}
for the basics of the jump; and \cite[Chapters 4-5]{AK00} for iterating the jump
through the computable ordinals, which generalizes by relativization to all countable ordinals.

When we prove it, Lusin's Theorem will reappear as Theorem \ref{thm:Lusin}.
It is normally stated as above, for the real numbers $\R$, but
it also applies to Cantor space $2^\N$.  Additionally, there are analogous theorems
(for both $\R$ and $2^\N$) using Baire category in place of Lebesgue measure.
Each of these is simpler to prove than Theorem \ref{thm:Lusin}, and we will use them
as warm-ups, introducing several of the techniques required in the more complex
context of Theorem \ref{thm:Lusin} itself.  Moreover, these results will demonstrate
the connection between computability and a standard question about these theorems.
The question is, ``why does Lusin's Theorem only yield $\g=\f$ up to a set of positive
measure, with different functions $\g$ for different positive values of $\ep$,
when the Baire-category version gives a single continuous $\g$ that equals $\f$
everywhere except on a meager set?''  This will be seen to be a direct consequence
of the nature of computing $X'$ from $\emptyset'\oplus X$ uniformly.  A uniform
computation is possible on a comeager set of $X$-values (with no wrong answers
even for the other $X$, but sometimes with no answer given).  Under Lebesgue
measure, some guessing is required, and therefore each uniform computation
will be wrong on some set of positive measure (although at least it will always
provide an answer, albeit occasionally an incorrect one).  In some respects
this can be seen to reflect the distinction between the modern notions of
\emph{generic computability} and \emph{coarse computability}, as described
for instance in \cite{JS12}.

Notation here is standard and follows Soare \cite{S87}.  In particular, capital Greek
letters such as $\Phi$, $\Gamma$, and $\Psi$ denote oracle Turing functionals,
and $\Phi_0,\Phi_1,\ldots$ is a standard computable enumeration of all such
functionals.  Sometimes $\Phi_{e,s}$ is used to denote the version of
the functional $\Phi_e$ that automatically halts after $s$ steps, even if it
has no output yet.  We write $a,b,q,r, u,v$ and sometimes $\ep$
for rational numbers; $\a,\b,\r,\x,\y$ for real numbers; capitals including
$A$, $C$, $X$, and $Y$ for elements of Cantor space $2^\N$ (that is,
for subsets of $\N$, often used as oracles for a functional), and $\f$, $\g$ for
functions mapping $2^\N\to 2^\N$ or $\R\to\R$.
Some further notation is described in the appendices
(Sections \ref{sec:computableanalysis} and \ref{sec:ordinaljumps}).

\section{Computing Discontinuous Functions}
\label{sec:main}

Section \ref{sec:computableanalysis} gives background in computable
analysis, including our reasons for using enumerations of Dedekind cuts
to name real numbers, rather than fast-converging Cauchy sequences.
Theorem \ref{thm:Weihrauch} shows that it is impossible to compute a discontinuous
function $\f:\R\to\R$ just from enumerations of Dedekind cuts for the input $\x\in\R$.
However, if we allow ourselves more information about $\x$, then it becomes possible.

 \begin{defn}
 \label{defn:jumpcomputable}
 A function $\f:2^\N\to 2^\N$  is \emph{jump-computable} if there
 exists a Turing functional $\Phi$ such that, for every $X\in 2^\N$,
 $\Phi^{(X')}:\N\to\{ 0,1\}$ is the characteristic function of $\f(X)$.
 Likewise, $\f:\R\to\R$ is \emph{jump-computable} if there
 exists a Turing functional $\Phi$ such that, for every $\x\in\R$
 and every enumeration $X=A\oplus B$ of the Dedekind cuts of $\x$, the function
 $$ \Phi^{(X')}:(\Q\times\N)^2\to\{ 0,1\}$$
computes an enumeration of the Dedekind cuts of $\f(\x)$.

More generally, for an ordinal $\alpha<\omega_1$ and an oracle set
$S\subseteq\N$, $f$ is \emph{$\alpha$-jump $S$-computable}
if there exists $\Phi$ such that, for every $\x\in\R$ and every enumeration
 $X$ of the Dedekind cuts of $\x$,
$$ \Phi^{(S\oplus X)^{(\alpha)}}:(\Q\times\N)^2\to\{ 0,1\}$$
computes an enumeration of the Dedekind cuts of $\f(\x)$;
similarly for functions on $2^\N$.
\end{defn}

If $\alpha$ is a countable noncomputable ordinal, then the $\alpha$-th
jump is not well-defined in general.  Section \ref{sec:ordinaljumps} explains
how we can choose an oracle $S$ complex enough to compute the complete diagram
of a presentation $\A$ of $\alpha$ (i.e., a linear order $\A$ isomorphic to
$(\alpha,\in)$ whose domain is $\N$).
So it is possible to discuss the situation $\alpha\geq\omega_1^{CK}$,
although one must fix an $S$-decidable presentation $\A$ of $\alpha$.
Since we defined the $\A$-jump $C^{(\A)}$ to have $C$ itself
as its $0$-th column, a functional $\Phi$ with oracle $(S\oplus X)^{(\A)}$
can recover $S$ from that column of the oracle, use it to compute
the presentation $\A$, and thus make sense of the oracle $(S\oplus X)^{(\A)}$ uniformly.

The principal theorem relevant here can be found in \cite{Kechris}.
By Lemma \ref{lemma:conversion}, for functions on $\R$,
this theorem holds with $\x$ and $\f(\x)$ represented either by fast-converging
Cauchy sequences or by enumerations of Dedekind cuts.
\begin{thm}
\label{thm:Borel}
For every Borel function $\f:2^\N\to 2^\N$, there is a Turing functional $\Phi$,
an oracle $S\subseteq\N$, and an $S$-decidable presentation $\A$
of some countable ordinal such that, for every $X\in 2^\N$,
$\Phi^{(S\oplus X)^{(\A)}}:\N\to\{ 0,1\}$ is the characteristic function of $\f(X)$.

Similarly, for every Borel function $\f:\R\to\R$, there is a Turing functional $\Phi$,
an oracle $S\subseteq\N$, and an $S$-decidable presentation $\A$
of some countable ordinal such that, for every enumeration $A\oplus B$ of the
Dedekind cuts of any $\x\in\R$, $\Phi^{(S\oplus X)^{(\A)}}$ enumerates
the Dedekind cuts of $\f(\x)$.
\end{thm}

Computable analysts have commonly approached jump-computability
the opposite way, by taking the limit of a computable function:
$$ \f(\x) = \lim_{s\to\infty}\text{$\Phi^X$}\!(~\cdot~, s).$$
That is, for each $s\in\N$, one computes an approximation to $\f(\x)$,
and the actual value $\f(\x)$ is the limit of these approximations.
As an example, the derivative of a differentiable function $\h(\x)$
could be given by letting $\Phi^X$ be the difference quotient
$s\cdot(\h(\x+\frac1{\text{$s$}})-\h(\x))$.  The Turing functional
$\Phi$ is readily defined from the functional and oracle computing $\h$.
(The derivative of a computable differentiable function is not in general
computable, so this is often the best that can be done.)
The connection between this method and ours is given by the
\emph{Transparency Lemma}.
\begin{lemma}[Transparency Lemma (folklore)]
\label{lemma:transparency}
A function is jump-computable if and only if it is the limit
of a computable function (in the sense immediately above).
\end{lemma}
So one could attack Lusin's Theorem by iterating the limit
operation instead of the jump operation.  Our choice to use the jump
is dictated mainly because it allows us to apply the known results
of Section \ref{sec:ordinaljumps} on near-uniform continuity of the jump,
especially Theorem \ref{thm:uniformity}.

\section{Results on Cantor Space $2^\N$}
\label{sec:Cantorspace}

Our two main theorems about Borel functions on Cantor space
depend heavily on three lemmas, which are the heart of the connection
between computability and Lusin's Theorem for $2^\N$.
These lemmas are well-known, but due to their importance, we give the proofs here,
at least for the base case $\alpha=1$.  Understanding them will prepare
the reader for the analogous theorems about functions on $\R$,
which are conceptually similar but more technical.  The three analogous
lemmas for $\R$ all appear in the appendix (Section \ref{sec:ordinaljumps}).

\begin{lemma}
\label{lemma:generic}
Let $X\in 2^\N$ be generic.  Then $X'\equiv_T \emptyset'\oplus X$,
and indeed there exists a single Turing functional $\Psi$ such that,
for every generic (or even $1$-generic) $X\in 2^\N$, $\Psi^{(\emptyset'\oplus X)}$
computes the characteristic function of $X'$.

More generally, for every oracle set $S\subseteq\N$ and every $S$-decidable
presentation $\A$ of a countable ordinal $\alpha$, there exists some $\Psi$
such that, for all $X\in 2^\N$ that are generic relative to $S$,
$\Psi^{(S^{(\A)}\oplus X)}$ computes the characteristic function of $(S\oplus X)^{(\alpha)}$.
Moreover, even for non-generic $X\in 2^\N$ (and for every $e\in\N$),
$\Psi^{(S^{(\A)}\oplus X)}(e)$ either diverges or computes
correctly whether $e\in (S\oplus X)^{(\alpha)}$.
\end{lemma}
\begin{pf}
We describe here the situation with $S=\emptyset$ and $\alpha=1$.
Given the oracle $(\emptyset'\oplus X)$, the functional $\Psi$ on input $e\in\N$
asks, for $s=1,2,\ldots$ in turn:
\begin{itemize}
\item
whether $\Phi_{e,s}^{(X~\!\res~\!\!s)}(e)$ halts (using no more of the oracle than the first $s$ bits); and also
\item
whether there exists any $\sigma\supseteq X\res s$ and any $t$ such that $\Phi_{e,t}^\sigma(e)$ halts.
\end{itemize}
If the first answer is ever positive, $\Psi$ concludes that $e\in X'$, while if the second answer is ever
negative, $\Psi$ concludes that $e\notin X'$.  Each conclusion is clearly correct (if ever reached),
and for a $1$-generic $X$, one of these conclusions must eventually be reached.
(Notice also that, even for non-generic $X$, $\Psi^{(\emptyset'\oplus X)}$ never gives
a wrong answer.  However, it could simply never give an answer.)
\qed\end{pf}

For successor $\alpha=\beta+1$, a similar
procedure relativized to $\emptyset^{(\beta)}$ succeeds.
The proof of the general result uses the uniformity of the above procedure
for all $\beta<\alpha$.

Our first Lusin-type result follows from Lemma \ref{lemma:generic}.
This is a standard theorem.
\begin{thm}
\label{thm:genericCantor}
Let $\f:2^\N\to 2^\N$ be a Borel function.  Then $\f$ restricts to a continuous function
$\g$ whose domain is the (comeager) set $G$ containing all elements of $2^\N$
that are generic relative to $S$.
\end{thm}
\begin{pf}
Theorem \ref{thm:Borel} yields $\Phi$, $S$, and $\A$ (as described there)
such that $\Phi^{(S\oplus X)^{(\A)}}$ computes $\f(X)$ for every $X\in 2^\N$.
Its restriction to $G$ is therefore computed by
$$ \Phi^{\Psi^{(S^{(\A)}\oplus X)}},$$
as witnessed by Lemma \ref{lemma:generic}.  
Since this function is $S^{(\A)}$-computable, it is continuous on its domain,
which contains the comeager set $G$.
(For $X\notin G$, either the computation will output the correct value $\f(X)$,
or there will exist some input $n$ on which the computation never halts.
As noted above, the computation will never give an incorrect answer, even for a single $n$.)
\qed\end{pf}

Now we turn to Lebesgue measure on $2^\N$, for which each
basic open set $\U_\sigma=\set{X\in 2^\N}{\sigma\subseteq X}$
has measure $2^{-|\sigma|}$, with the measure extending to all measurable
subsets of $2^\N$ in the usual Lebesgue definition.  There is
an analogue of Lemma \ref{lemma:generic} for Lebesgue measure,
but, perhaps contrary to one's expectations, it is uniform only for
sets of arbitrarily small measure, not up to a set of measure zero.
The following lemma appears as Theorem 2 of \cite{S72}, but a proof
may aso be deduced from the (more complicated) Lemma \ref{lemma:GL1}
and Theorem \ref{thm:uniformity} in the appendices, especially in concert
with the proof of Lemma \ref{lemma:errorset} below.

\begin{lemma}
\label{lemma:LebesgueCantor}
For every oracle set $S\subseteq\N$ and every $S$-decidable
presentation $\A$ of a countable ordinal $\alpha$, there exists some $\Psi$
such that, for each fixed rational $\ep>0$, 
$\Psi^{(S^{(\A)}\oplus X)}(\ep,~\cdot~)$ computes the characteristic function of
$(S\oplus X)^{(\A)}$
%the statement
%$$(\forall e\in\N)~[\chi_{(S\oplus X)^{(\A)}}(e)=\Psi^{(S^{(\A)}\oplus X)}(\ep,e)]$$
%holds 
for all $X$ outside a subset of $2^\N$ of measure $<\ep$.
Moreover, the computation $\Psi^{(S^{(\A)}\oplus X)}(\ep,e)$ %on the right side 
halts for every $X\in 2^\N$ and every $(\ep,e)$,
although for $X$ within the set of measure $<\ep$ it may output
a value distinct from the ``correct'' answer $\chi_{(S\oplus X)^{(\A)}}(e)$.
\end{lemma}
We will write $\Psi_\ep$ for the function $\Psi(\ep,~\cdot~)$, with any oracle.
For Lebesgue measure, a second fact is also necessary.

\begin{lemma}
\label{lemma:errorset}
In Lemma \ref{lemma:LebesgueCantor}, the error set described, whose measure
is $<\ep$, is the union of an $S$-computably enumerable sequence of basic open
sets, uniformly in both $\A$ and $\ep$.
\end{lemma}
\begin{pf}
First we explain this proof when $\alpha=1$ and with $S$ omitted.  Now
$\Psi_\ep^{\emptyset'\oplus X}$ is attempting to compute $X'$ as described
in Lemma \ref{lemma:LebesgueCantor}.  Whenever
$\Psi_\ep^{\emptyset'\oplus X}(e)\converges=1$, it does so because it has seen
the computation $\Phi_e^X(e)$ converge already, so $e\in X'$ and this answer
is correct.  All errors (if any occur for this $X$) occur for values $e$ where
$\Psi_\ep^{\emptyset'\oplus X}(e)\converges=0$.

To enumerate an error set, therefore, once $\Psi_\ep$ has declared that
$\Psi_\ep^{\emptyset'\oplus X}(e)\converges=0$ for all $X$ extending some
particular $\tau\in 2^{<\N}$, we watch for any strings $\rho\supseteq\tau$
and any $t$ for which $\Phi_e^{\rho}(e)\converges$.  Any time we find such a $\rho$,
we enumerate the basic open set defined by $\tau$ into our error set, knowing that
this entire basic open set lies in the error set.  Since all errors are errors of this type
(guessing that $e\notin X'$ and later being proven wrong), this list of basic open sets
is precisely the error set and thus has measure $<\ep$.

In the general setting, the jump $(S\oplus X)^{(\A)}$ consists of many columns,
one for each point $k$ in the presentation $\A$ of the ordinal $\alpha$.  If $k$
represents a successor ordinal $\beta+1$ in $\A$, then we do the same procedure as above,
waiting for $\Psi_\ep$ to use its own computation of $(S\oplus X)^{(\beta)}$
to declare that $e\notin (S\oplus X)^{(\beta+1)}$ for all $X\supset\tau$ and then watching
for $\Phi_e^{(S\oplus \rho)^{(\beta)}}(e)\converges$ for each $\rho\supseteq\tau$.
Only $\ep\cdot 2^{-k-1}$ of the measure $\ep$ of the error set is alloted to this process,
and it certainly enumerates all $X$ for which $\Psi_\ep$ gave the wrong answer about
whether some $e$ lies in $(S\oplus X)^{(\beta)}$ but gave correct answers to these
questions for all ordinals $<\beta$.  Now if in fact $\Psi_\ep$ gives a wrong answer for some
$e$ about whether $e\in (S\oplus X)^{(\A)}$, then there is a least element $\beta$
of $\alpha$ for which it did so, and for this $\beta$, $X$ will belng to a basic open set
enumerated into our error set.  
\qed\end{pf}

\begin{thm}
\label{thm:LebesgueCantor}
Let $\f:2^\N\to 2^\N$ be a Borel function, and fix $\ep>0$.
Then there exists a continuous function $\g:2^\N\to 2^\N$
such that 
$$  \mu ( \set{X\in 2^\N}{\g(X)\neq \f(X)} ) <\ep.$$
\end{thm}
\begin{pf}
The proof is somewhat analogous to that of Theorem \ref{thm:genericCantor},
using Lemmas \ref{lemma:LebesgueCantor} and \ref{lemma:errorset}.  In order to ensure
that $\g$ is defined on all of $2^\N$ (which was impossible in
Theorem \ref{thm:genericCantor}), we must provide a value for $\g(X)$
even for those $X$ such that $\Psi_\ep^{(S^{(\A)}\oplus X)}\neq (S\oplus X)^{(\A)}$.
The danger here is that for these $X$, the functional $\Phi$ running with the
incorrect oracle $\Psi_\ep^{(S^{(\A)}\oplus X)}(\ep,~\cdot~)$ may fail to compute
a total function, in which case our naive version of $\g$, namely
$$ \Phi^{\Psi_\ep^{(S^{(\A)}\oplus X)}},$$
will fail to map all of $2^\N$ into $2^\N$.
Therefore it is necessary for us to be able to realize, at some point during the computation
above, that $\Psi_\ep^{S^{(\A)}\oplus X}$ fails to compute $(S\oplus X)^{(\A)}$ correctly.
This is exactly the point of Lemma \ref{lemma:errorset}.

The procedure for $\Gamma$ begins by starting the naive computation
$\Phi^{\Psi_\ep^{(S^{(\A)}\oplus X)}}$.  If we ever reach a stage $s$ at which
the $s$-th basic open set (given by a string $\sigma_s\in 2^{<\N}$
such that the basic open set is $\set{Y\in 2^\N}{\sigma\subset Y}$)
in the enumeration above contains $X$ -- i.e., $\sigma_s\subset X$ --
then at that stage we simply end the computation on all inputs:
the output is
$$ \g(X) = \set{n\in\N}{\Psi_{\ep,s}^{(S^{(\A)}\oplus X)}(n)\converges=1},$$
where $\Psi_{\ep,s}$ indicates that we only run $\Psi_\ep$ for $s$ steps.
If we never reach such a stage $s$, then $\Gamma$ continues
running each computation $\Psi_\ep^{(S^{(\A)}\oplus X)}$ until it halts,
as it must, since in this case, with $X$ not in the error set, this function
succeeds in computing $\f(X)$.  Thus $\Gamma$ does indeed
compute a function $\g$ defined on all of $2^\N$, and thus
continuous, differing from $\f$ only on the error set.  Notice also
that the computation of an index for $\g$ is uniform in the positive rational $\ep$ chosen.
\qed\end{pf}

Theorems \ref{thm:genericCantor} and \ref{thm:LebesgueCantor}
are not perfectly analogous.  The former yields a single continuous version of $\f$
that may be undefined on a small (i.e., meager) set of inputs $X$.
The latter yields different versions $\g_{\ep}$ for each $\ep>0$,
each of which is continuous on all of $2^\N$ but may differ
from $\f$ on a small set (i.e., of measure $<\ep$).  The computability-theoretic difference
between working with generic sets and working with Lebesgue measure
accounts for the difference in the corresponding lemmas, and thus may be seen
as explaining the difference between the theorems.

The technique of enumerating the error set, used in Lemma \ref{lemma:errorset}, will reappear in Section
\ref{sec:Lusin}, when we address Lusin's Theorem itself using Lebesgue measure on $\R$.

\section{Baire Category on $\R$}
\label{sec:Bairecat}

Our next step is to address the real line $\R$ rather than Cantor space $2^\N$.
The difficulty here lies in representing real numbers $\x\in\R$, since our fundamental
objects are elements $X\in 2^\N$ and each $\x\in\R$ will be represented by many different
$X\in 2^\N$, each naming a Cauchy sequence in $\Q$ that converges fast to $\x$.
(In Section \ref{sec:Cantorspace}, this was not an issue:
each $X\in 2^\N$ simply represented itself!)  Indeed, the $X$'s representing
a single $\x$ will have many distinct Turing degrees.

The first impulse, for addressing this problem, is to switch to Dedekind cuts,
since each $\x\in\R$ is represented by a single strict left Dedekind cut (i.e.,
a nonempty subset of $\Q$, downward-closed under $<$, with no greatest element).
Unfortunately, these cannot be used directly as oracles:  the analogue
of Lemma \ref{lemma:conversion} for this system of representation of real numbers
is false.  Lemma \ref{lemma:conversion} is our compromise, converting
fast-converging Cauchy sequences into enumerations of strict left and right Dedekind cuts
(and back) effectively.  This is not a perfect solution, since for a given $\x$,
the strict left and right cuts $L_\x$ and $R_\x$ each have many distinct enumerations,
of many Turing degrees.  However, as we now describe, it enables us to use
enumerations of these cuts to run oracle Turing computations in an effective way,
so that the output depends only on $L_\x$ and $R_\x$ (and thus only on $\x$ itself),
not on the choice of enumerations of $L_\x$ and $R_\x$.
Definition \ref{defn:enumcut} describes the ``canonical'' enumeration of each $\x$,
which is crucial for the proof of the next theorem.  Our use of Dedekind cuts
stems from our lack of any similar notion of a canonical Cauchy sequence converging fast to $\x$.

The generic real numbers form a comeager subset $G$ of $\R$, where by definition
$\x\in\R$ is generic just if the cut $L_\x$ is generic among all downward-closed
sets of rational numbers. Clearly, genericity among such sets rules out having a greatest
element, being empty, or being co-empty, so these sets are all strict left cuts of irrational numbers $\x$.
Moreover, no definable left cut $\x$ can be generic in this sense.
This set $G$ will be the comeager set we use to prove the Baire-category
version of Lusin's Theorem.

\begin{thm}
\label{thm:genericR}
Let $\f:\R\to\R$ be a Borel function.  Then $\f$ restricts to a continuous function
$\g$ whose domain is the (comeager) set $G$ containing all elements of $\R$
that are generic relative to $S$.
\end{thm}
\begin{pf}
Once again, Theorem \ref{thm:Borel} yields $\Phi$, $S$, and $\A$
such that, for every $\x\in\R$ and every enumeration $A\oplus B$ of the
Dedekind cuts of $\x$,
$\Phi^{(S\oplus A\oplus B)^{(\A)}}$ enumerates the Dedekind cuts of $\f(\x)$.
On an enumeration $A\oplus B$ of $\x\in G$, one might hope for it to be computed
by a function along the lines of
$$ \Phi^{\Psi^{(S^{(\A)}\oplus (A\oplus B))}},$$
as witnessed by Lemma \ref{lemma:generic}.  
Since this function is $S^{(\A)}$-computable, it would be continuous on its domain.

Unfortunately, this does not suffice for a proof.  This function is only defined on generic
enumerations of cuts, and no matter which $\x\in\R$ one chooses, there will exist
non-generic enumerations of its cuts.  Thus, our program above may fail to accept
certain enumerations of $\x$, whereas a computable function must accept all of them.
To rectify this, we give a slightly more involved program, which will succeed
in computing $\f(\x)$ below an $S^{(\A)}$ oracle for every enumeration of each $\x\in\R$
whose left Dedekind cut $L_\x$ forms an $S$-generic downward-closed subset of $\Q$.
We first describe our program $\Gamma$, which ``uniformizes'' the enumerations
$A\oplus B$ given to it, and then analyze it.

$\Gamma$ is given the oracle $A\oplus B$, along with the fixed set $S^{(\A)}$.
At each stage $s$, it defines
$$ A_s=A\cap(\{0,\ldots,s\})^2,~B_s=B\cap(\{0,\ldots,s\})^2,\text{~and~}
l_s = \min \set{n\in\N}{q_n\notin \pi_1(A_s\cup B_s)},$$
the greatest (up to length $s$) initial segment of a fixed computable listing
$\{ q_0,q_1,\ldots\}$ of $\Q$ for which the cuts $L_\x$ and $R_\x$ are fully
defined by stage $s$.  Then it runs the program $\Phi$ for $s$ steps
on the given input $q\in\Q$, using the oracle:
$$ \Psi_s^{S^{(\A)}\oplus (\pi_1(A_s)\times\N)\oplus (\pi_1(B_s)\times\N)}$$
and outputs the same value as this program
in case it halts, or else goes back and repeats the same procedure for $s+1$.
(If, in response to an oracle question from $\Gamma$, $\Psi$ fails to halt within $s$ steps
using its own smaller oracle, then again $\Gamma$ starts its own program over at $s+1$.)
This is the entire program for $\Gamma$.

In $s$ steps, $\Psi$ can never ask an oracle question about an element $>s$,
so all oracle values used by $\Psi$ here are exactly as in the canonical enumeration
$(L_\x\times\N)\oplus(R_\x\times\N)$ of the cuts of $\x$.  Thus, if $\Psi$
actually gives an answer (on an input $e$) using its oracle here, that answer
will be identical to the one given by $\Psi$ using that canonical oracle.
In case $\x$ is generic relative to $S$, it will therefore tell whether
$e\in (S\oplus (L_\x\times\N)\oplus (R_\x\times\N))^{(\A)}$
or not.  Thus our uniformization to the canonical oracle will have eliminated
all uncertainty introduced by the use of the enumeration $(A\oplus B)$.
Moreover, if $l_s$ goes to $+\infty$ as $s$ increases, then $\Psi$ will ultimately
have access to as much of the canonical oracle as it needs for its computation.
Finally, notice that a generic $\x\in\R$ cannot be a rational number, and therefore
every $q_n\in\Q$ lies in $L_\x\cup R_\x$, hence lies in $\pi_1(A_s\cup B_s)$
for all sufficiently large $s$.  Thus for generic $\x$, we automatically have
$\lim_s l_s=+\infty$, and so $\Psi$ ultimately computes (correctly) as much
of $(S\oplus (L_\x\times\N)\oplus(R_\x\times\N))^{(\A)}$ as $\Phi$
requires.  In short, for every generic $\x\in\R$ \emph{and every enumeration $A\oplus B$
of the cuts of $\x$}, the program $\Gamma$ produces the exact same output as 
$\Phi^{(S\oplus A\oplus B)^{(\A)}}$, namely, an enumeration of the left
and right cuts of $\f(\x)$, exactly as claimed by the theorem.

Finally, notice that $\Psi$ always uses an oracle that describes
$(L_\x\times\N)\oplus(R_\x\times\N)$ correctly:  the only possible problem
arises when $\x$ is not $S$-generic, and even then, $\Psi$ never gives a wrong answer.
Therefore, $\Gamma$ never outputs a wrong answer either:  for each $q\in\Q$ and
each enumeration $A\oplus B$ of any $\x\in\R$, if $\Gamma$ halts and outputs
$0$ on input $q$, then indeed $q\in L_{\f(\x)}$; likewise the output $1$ always
indicates that $q\in R_{\f(\x)}$.  For an $\x$ that is not generic relative to $S$, the program may
fail to halt on various inputs $q\in\Q$, in which case that $\x$ is not in the domain
of the continuous function defined by $\Gamma$.  However, non-$S$-generic real numbers
form a meager subset of $\R$, so we have proven the theorem.
\qed\end{pf}

\section{Lusin's Theorem}
\label{sec:Lusin}

Now we may approach Lusin's Theorem.  The key lemma,
again well-known, is described in more
detail in the second appendix (Section \ref{sec:ordinaljumps}),
as Corollary \ref{cor:uniformity}.

\begin{lemma}
There exists a Turing functional $\Psi_{\ep}$, uniform in the rational $\ep>0$,
such that for every set $S\subseteq\N$ and every fixed $S$-decidable
presentation $\A$ of any ordinal $\alpha<\omega_1^S$
(with $\alpha$-jumps $C^{(\A)}$ defined using this presentation),
the ``error set''
$$ \U_{\ep,S, \A}=\set{\x\text{$~\!\in\R$}}{\text{$\Psi_\ep^{S^{(\A)}\oplus L_{\x}\oplus R_{\x}}\neq
(S\oplus L_{\x}\oplus R_{\x})^{(\A)}$}}$$
has measure $<\ep$ and is an $S^{(\A)}$-effective union of rational open intervals,
uniformly in $\ep$, $S$, and $\A$.
Moreover, for all $\x$, $\Psi_\ep^{S^{(\A)}\oplus L_{\x}\oplus R_{\x}}$ is total.
\qed\end{lemma}

\begin{thm}[Lusin, 1912]
\label{thm:Lusin}
Every Borel-measurable function $\f:\R\to\R$ is nearly continuous.
\end{thm}

\begin{pf}
As usual, Theorem \ref{thm:Borel} provides $\Phi$, $S$, and $\A$
such that, for every $\x\in\R$ and every enumeration $A\oplus B$ of the
Dedekind cuts of $\x$,
$\Phi^{(S\oplus A\oplus B)^{(\A)}}$ enumerates the Dedekind cuts of $\f(\x)$.
We wish to construct a $\Gamma$, similar to that of Theorem \ref{thm:genericR},
that uses $\Psi$ to compute $(S\oplus L_\x\oplus R_\x)^{(\A)}$ from $(S^{(\A)}\oplus A\oplus B)$,
then applies $\Phi$ to enumerate the cuts of $\f(\x)$.  The catch is that now
we wish the $\g$ computed by $\Gamma$ to have domain $\R$; on the other hand,
we may allow $\g(\x)\neq\f(\x)$ on the error set $\U_{\ep,S,\A}$ given by
Corollary \ref{cor:uniformity}.  We write $\U=\U_{\ep,S,\A}$, having fixed the
three parameters, and, using Corollary \ref{cor:uniformity},
fix an $S^{(\A)}$-computable enumeration of rational open intervals
whose union equals $\U$.  We also fix a computable list $q_0,q_1,\ldots$
of all rational numbers, without repetitions.  It will be convenient to assume
that the interval of $\U$ enumerated at stage $s$ is of the form $(q_m,q_n)$
with $m,n\leq s$.  This is not difficult to arrange, except that we must allow
there to be stages at which no interval is enumerated.

The crucial fact here is the $S^{(\A)}$-computable enumerability of
the error set $\U$ as a union of intervals with rational end points.
This will allow our $\Gamma$ to recognize these intervals as they appear,
at which stage it will abandon its hope of computing a function equal to $\f$ within such an interval,
and will instead do ``damage control'' to ensure that it does actually compute
a function (necessarily continuous) on each such interval.  As long as we
make these ``damage control'' functions meet at the end points of their intervals,
we will have continuity everywhere:  outside the error set, $\Gamma$ will compute
$\f(\x)$ successfully, in exactly the style of Theorem \ref{thm:genericR},
since $\Psi$ makes no mistakes outside the error set in computing
$(S\oplus (L_\x\times\N)\oplus (R_\x\times\N))^{(\A)}$.

It complicates matters, but it is also liberating, to recall that every rational $\x$
must lie in $\U$.  (In particular, $\Psi$ will always guess incorrectly
about whether the program index $e$ lies in $(S\oplus A\oplus B)^{(\A)}$, where, for every input,
$\Phi_e$ halts as soon as the rational $\x$ appears in $\pi_1(A\cup B)$.)
The complication lies in the fact that, whenever a rational interval $(a,b)$
is enumerated into the error set, the rational $a$ itself must lie in another
interval of the error set, as must $b$, causing these intervals ultimately to metastasize
into one long open interval with irrational endpoints, the union of infinitely many
overlapping rational intervals in $\U$.  The liberation results from the
fact that, as in Theorem \ref{thm:genericR}, we are again free to ignore the difficulty
that, for rational $\x$, $L_\x\cup R_\x$ omits an element (namely $\x$ itself).  Given an
enumeration $A\oplus B$ of the cuts of an arbitrary $\x\in\R$, we know that either
$\x$ will eventually enter the error set (which we can enumerate!), or else
$L_\x\cup R_\x=\Q$.  As long as $\x$ is not yet in the error set, therefore, we are
safe in waiting arbitrarily long for the next rational $q_n$ on our list to be enumerated
by either $A$ or $B$.  If $\x$ ever does appear in the error set, we switch
to damage control on $\x$.

The bulk of our description of the procedure for the functional $\Gamma$
deals with the rational intervals (finitely many of them, at each stage $s$)
that currently constitute the error set.  Once we have settled what to do with them,
the procedure on the remainder of $\R$ is the same one already applied several times.
So suppose that as of stage $s$, our enumeration of $\U$ so far consists of
rational intervals $(a_1,b_1),~(a_2,b_2),\ldots,(a_{s-1},b_{s-1}$, with each $b_i < a_{i+1}$.
By induction we assume that we have \emph{assigned values} (as defined below)
to all points in the closure $\overline{\U_s}$.  This means that we have already determined the
finitely many values $\g(a_i)$ and $\g(b_i)$, possibly along with values $\g(c)$
at finitely many other rational $c\in\U_s$, and that we intend $\g$ on $\overline{\U_s}$
to be the piecewise linear function that ``connects these dots'' within each connected
component of $\overline{\U_s}$.

Before we describe the action of $\Psi$ on a particular $(A\oplus B)$, we explain how it
sets the scene at each stage $s+1$.  We have a fixed computable enumeration
$\{ q_0,q_1,\ldots\}$ of $\Q$, so suppose that $\{ q_0,\ldots,q_s\}$ is ordered as
$q_{j_0}<\cdots<q_{j_s}$.  These define $s+1$ open intervals $(-\infty,q_{i_0})$,
$(q_{j_0},q_{j_1})$, etc., each giving an initial piece of a Dedekind cut:
write $L_{i,s}=\{ q_{j_0},\ldots,q_{j_{i-1}}\}$ and $R_{i,s}=\{ q_{j_i},\ldots,q_{j_{s}}\}$.
$\Psi$ runs the programs
$$\Phi_s^{S^{(\A)}\oplus (L_{i,s}\times\N)\oplus (R_{i,s}\times\N)}(q_k)$$
for each $i\leq s$ and each $k\leq s$, giving up after $s$ steps if there is not yet
any convergence.  This gives us general constraints on the possible values of $\f$ in each
interval $(q_{j_{i-1}},q_{j_i})$:  an upper bound (possibly $+\infty$)
$$ u_{i,s} = \min\set{q_k}{k\leq s~\&~\Phi_s^{S^{(\A)}\oplus (L_{i,s}\times\N)\oplus (R_{i,s}\times\N)}(q_k)\converges=1}$$
and a lower bound (possibly $-\infty$)
$$ v_{i,s} = \max\set{q_k}{k\leq s~\&~\Phi_s^{S^{(\A)}\oplus (L_{i,s}\times\N)\oplus (R_{i,s}\times\N)}(q_k)\converges=0}.$$
If the program in question actually computes $\f$ on $(q_{j_{i-1}},q_{j_i})$ correctly,
then the values of $\f$ there must all lie in $(v_{i,s},u_{i,s})$.  Thus the set of boxes
$(q_{j_{i-1}},q_{j_i}) \times (v_{i,s},u_{i,s})$ may be viewed as a tentative approximation
of the graph of $\f$ in the $xy$-plane.  Of course, this does not take the error set $\U$
into account.  It could happen that $u_{i,s}\leq v_{i,s}$ under this definition, if the interval
$(q_{j_{i-1}},q_{j_i})$ is contained within $\U$.  Also, two adjacent boxes (as described above)
might have nonintersecting closures, which again would indicate interference by the error set.
Therefore we have only defined $u_{i,s}$ and $v_{i,s}$ here, without using them as yet.

We now continue stage $s+1$ by considering
the next interval $I_s=(a_s,b_s)$ to appear in the enumeration
of intervals comprising $\U$.  We wish to assign $\g$-values to all
points $\x\in [a_s,b_a]$.  This is the damage-control operation:
we have given up on computing $\f(\x)$ correctly for these points,
and merely hope to make $\g$ continuous.  Within $I_s$ this will be
easy, but we also need to ensure that our chosen values allow $\g$
to be continuous at the boundary of $\U$, where it interfaces with
error-free computations of $\f$.  (The boundary of $\U$ is in fact
all the rest of $\R$, since $\U$ is dense.)

Let $\U_{s+1}=\cup_{i<s} [a_i,b_i]$. If $I_s\subseteq\U_s$
(which is decidable from the end points), we do nothing, since values
will have already been assigned to all points in $\U_s$ at preceding stages.
Otherwise $(I_s-\U_s)$ is a finite union of closed nonempty intervals,
and we treat each of its connected components $[c,d]$ separately.
First, if any of the points $q_{j_0},\ldots,q_{j_s}$ defined earlier lies in $[c,d]$,
then for each such $q_{j_i}$ satisfying $\min(u_{i,s},u_{i,s+1})\geq \max(v_{i,s},v_{i+1,s})$,
we define
$$ \g(q_{j_i}) = \frac{\min(u_{i,s},u_{i+1,s})+\max(v_{i,s},v_{i+1,s})}2.$$
(In case $\min(u_{i,s},u_{i+1,s})=+\infty$, $\g(q_{j_i})=\max(v_{i,s},v_{i+1,s})+1$.
In case $\max(v_{i,s},v_{i+1,s})=-\infty$, $\g(q_{j_i})=\min(u_{i,s},u_{i+1,s})-1$.
If both occur, then $\g(q_{j_i})=0.$)

If $q_{j_i}$ lies in $[c,d]$ but $\min(u_{i,s},u_{i,s+1}) < \max(v_{i,s},v_{i+1,s})$,
then these bounds create an inconsistency at $q_{j_i}$.  In this situation
we continue to enumerate $\U$ until a rational interval $(r_0,r_1)$ containing
$q_{j_i}$ is enumerated into $\U$.  (This must occur eventually, as all rational
numbers lie in $\U$.)  When this first happens, we define
$$ \g\left(\frac{\max(r_0,c,q_{j_{i-1}})+q_{j_i}}2\right) =\frac{u_{i,s}+v_{i,s}}2\hspace{10mm}
\g\left(\frac{q_{j_i}+\min(r_1,d,q_{j_{i+1}})}2\right) =\frac{u_{i+1,s}+v_{i+1,s}}2.
$$
The argument on the left lies below $q_{j_i}$ but above each of $r_0$, $c$, and $q_{j_{i-1}}$,
hence within $[r_0,r_1]$, within the interval where $u_{i,s}$ and $v_{i,s}$
are the upper and lower bounds, and within $[c,d]$, so it is safe to assign to it
this $\g$-value between those bounds; similarly for the argument on the right.
(In case $q_{j_i}=c\in\overline{\U_s}$, only the equation on the right applies,
as in this case the left equation has argument $c$ and $\g(c)$ is already defined.
Similarly, if $q_{j_i}=d\in\overline{\U_s}$, only the equation on the left applies.)

Next, we attend to the end points.
If either $c$ or $d$ lies in the closure $\overline{\U_s}$, then it already
has been assigned a $\g$-value.  If $c\notin\overline{\U_s}$, then
find the interval $(q_{j_{i-1}},q_{j_i})$ to which it belongs and define
$$\g(c) = \left\{\begin{array}{cl} \frac{u_{i,s}+v_{i,s}}2,&\text{~if this is defined;}\\
u_{i,s}-1,&\text{~if~}v_{i,s}=-\infty~\&~u_{i,s}\neq +\infty;\\
v_{i,s}+1,&\text{~if~}u_{i,s}=+\infty~\&~u_{i,s}\neq -\infty;\\
0,&\text{~if~}v_{i,s}=-\infty~\&~u_{i,s} = +\infty.
\end{array}\right.$$
Likewise if $d\notin\overline{\U_s}$
and $d\in (q_{j_{k-1}},q_{j_k})$, then define $\g(d)$ exactly the same way,
using $u_{k,s}$ and $v_{k,s}$.
The case where $c$ or $d$ actually equals $q_{j_i}$ for some $i\leq s$
was already covered in our first step.

Now we have defined $\g$ on $c$, on $d$, and on every $q_{j_i}$ in between.
We finish by ``connecting the dots'':  simply define $\g$ between each pair
of consecutive points among these by linear interpolation between the $\g$-values
of the two points.

It now remains to describe the procedure followed by $\Gamma$ with a specific
enumeration $A\oplus B$ of strict Dedekind cuts (of an arbitrary $\x\in\R$)
in its oracle.  We use sets $A_s$ and $B_s$ and the length $l_s$ just
as defined in the proof of Theorem \ref{thm:genericR}.
At stage $s+1$, $\Gamma$ searches for the least $t\in\N$
such that either
\begin{enumerate}
\item
$t > s$ and $l_t\geq s$; or
\item
there exist $m,n\leq t+1$
such that $q_m\in A_{t+1}$ and $q_n\in B_{t+1}$ and the open interval
$(q_m,q_n)$ is enumerated into $\U$ by stage $t+1$.
\end{enumerate}
If condition (2) holds, then by running the instructions above all the way up
to the stage $t$ at which $(q_m,q_n)$ is enumerated into $\U_t$,
we can determine the $\g$-values assigned at that stage to points in this
interval.  Of course, at this stage we still only have a finite approximation to $\x$,
but, knowing the finitely many linear pieces of $\g$ on this interval,
we can use the approximation given by $A_s$ and $B_s$ to approximate
$\g(\x)$ in the obvious way.  (Notice that this works even if $\x$ itself
is rational:  the fact that some single rational never appears in either Dedekind cut
does not stop us from computing $\g(\x)$.)

If condition (1) holds, then from $A_t$ and $B_t$ we can decide
which of the intervals $(q_{j_{i-1}},q_{j_i})$ defined at stage $s+1$ contains $\x$.  
If $u_{i,s}\leq v_{i,s}$, then we do nothing at this stage,
since then this entire interval must eventually enter $\U$.
Assuming $u_{i,s}<v_{i,s}$, we know that $\x$ itself does not lie in $\U_{s+1}$:
condition (1) requires $t>s$, so if an interval $(q_m,q_n)$ containing $\x$ had
been enumerated into $\U_{s+1}$, it would have $m,n\leq s$ (by our convention
for this enumeration), and thus stage $s$ would have satisfied condition (2).
However, we do not want to deal with an $\x$ that is too close to either
$q_m$ or $q_n$.  So we search again, for a stage $t'$ such that
$\U_{t'}$ contains both an interval $(a'_m,b'_m)$ that contains $q_m$
and also an interval $(a'_n,b'_n)$ that contains $q_n$, and such that either
\begin{enumerate}
\item[(1a)]
$A_{t'}$ contains a rational number $\geq b'_m$, and
$B_{t'}$ contains a rational number $\leq a'_n$
(making $\x\in (b'_m,a'_n)$); or
\item[(1b)]
$A_{t'}$ contains a rational number $\geq a'_n$ or
$B_{t'}$ contains a rational number $\leq b'_m$
(making $\x\in\U_{t'}$).
\end{enumerate}
If (1b) holds, $\Gamma$ does not converge on any inputs here,
because we know that our $\x$ will eventually enter $\U$.
If (1a) holds, then $\Gamma$ halts and outputs $1$ for all inputs $\geq u_{i,s}$,
and halts and outputs $0$ for all inputs $\leq v_{i,s}$.  (Recall that we
have checked that $v_{i,s}<u_{i,s}$.)  This means that these values
will be upper and lower bounds (respectively) for the real number $\g(\x)$
ultimately determined by $\Gamma$.  This completes our description
of the program executed by $\Gamma$.

It remains to show that $\Gamma^{S^{(\A)}\oplus A\oplus B}$ succeeds
in computing a function $\g:\R\to\R$, as $A\oplus B$ varies over all
enumerations of Dedekind cuts, and then that $\g(\x)=\f(\x)$ outside
of the error set $\U$ (which we know has measure $<\ep$).  Being computed
by $\Gamma$ this way, $\g$ will necessarily be continuous.

To see this, fix any $\x\in\R$ and any enumeration $A\oplus B$ of the cuts
$L_\x$ and $R_\x$ of $\x$.  Suppose first that $\x\in\U$.  Here we do not claim that $\g(\x)=\f(\x)$,
but we do show that $\Gamma^{S^{(\A)}\oplus A\oplus B}$ computes the same
real-number output $\g(\x)$ regardless of the enumeration of $\x$.
Indeed, since $\x\in\U$, there is some interval $I_{s_0}=(a_{s_0},b_{s_0})$
with $\x\in I_{s_0}$ enumerated into $\U$ at some stage $s_0$.
Notice that at every subsequent stage $s>s_0$, condition (2) will hold with $t=s_0-1$,
while condition (1) can only be applied at stage $s+1$ if $t>s$;
thus the instructions for condition (1) will never be followed again.
Moreover, regardless of the choice of enumeration $A\oplus B$
of $L_\x$ and $R_\x$, we will recognize from stage $s_0+1$ onwards
that the enumeration of $I_{s_0}$ into $\U$ has placed us in condition (2).
Therefore, from stage $s_0+1$ on, $\Gamma$ will always proceed
with the computation of the assigned value $\g(\x)$, based on the
linear functions chosen at stage $s+1$ for the interval of $\U_s$ containing $\x$.
Of course, the computation of these linear functions (with rational coefficients, no less!)
is effective regardless of the enumeration $A\oplus B$, and the
end points of the linear function were selected specifically so that no
upper or lower bounds enumerated by $\Gamma$ at preceding stages
could possibly contradict the value $\g(\x)$ assigned to any $\x$ in
the interval, including our specific $\x$.  (It is also important to note here
that we never assigned contradictory information in the computation
of $\g(\x)$, even when using condition (1).  All rationals enumerated into
the lower cut of $\g(\x)$ were checked to be $<$ all rationals enumerated
into the upper cut.  Moreover, as long as condition (1) applied, the same
rationals were enumerated into each cut independently of the specific
enumeration $A\oplus B$.  The choice of $A\oplus B$ only affects the stage
at which the computation permanently switches over to condition (2),
which occurs at or before stage $s_0$, and the scene-setting procedure
included a check that all upper and lower bounds enumerated for values
$\x$ in $I_{s_0}$ are consistent with the choice of $\g(\x)$.
Therefore $\Gamma$ does compute the same value $\g(\x)$ independent
of the specific enumeration $A\oplus B$ of the cuts of $\x$.  This is all
that we require of $\Gamma$ on an $\x\in\U$:  the damage-control
procedure succeeded.

The other case occurs when $\x\notin\U$.  Now the procedure above
never acts on condition (2).  Since every rational number lies in $\U$, it follows
that $\x\notin\Q$, so every $q_m$ lies in $L_\x\cup R_\x$, and thus the length
$l_s$ of the approximation $A_s\oplus B_s$ goes to $+\infty$ as $s$ increases.
Consequently, for every stage $s$, there does exist some $t>s$
with $l_t\geq s$.  

Now since $\x\notin\U$, $\Phi$ does enumerate the cuts of $\f(\x)$ when it runs
with the oracle set $\Psi^{S^{(\A)}\oplus (L_\x\times\N)\oplus (R_\x\times\N)}$.
Therefore, for any specific rational $u>\f(\x)$ and $v<\f(\x)$, there exists some $s_0$
such that $\Phi$ with this oracle enumerates $u$ into the upper cut and $v$
into the lower cut, in at most $s_0$ steps, using only the initial portion
$((L_\x~\!\res~\!s_0)\times\N)\oplus((R_\x~\!\res~\!s_0)\times\N)$
of the oracle of $\Psi$.  But $l_t\to +\infty$ as $t$ increases,
so at every stage $\geq s_0$, running $\Gamma$ with the oracle
$A\oplus B$ will produce this much of $L_\x$ and $R_\x$.
Now fix some $m,n>s_0$ such that $\x\in (q_m,q_n)$
but no $q_k$ with $k<\max(m,n)$ lies in $(q_m,q_n)$.  When we reach
stage $s+1$ with $s=\max(m,n)$, this interval will be the $i$-th interval
in the partition of $\R$ (that is, $m=j_{i-1}$ and $n=j_i$ at this stage),
and so the bounds $u_{i,s}$ and $v_{i,s}$ determined at that stage
will have $\f(\x) < u_{i,s}\leq u$ and $\f(\x)> v_{i,s}\geq v$.
Moreover, since $\x\notin\U$, condition (1b) must hold
(at this and all other stages), and so $u_{i,s}$ and $v_{i,s}$
are enumerated into the upper and lower cuts of $\g(x)$
at that stage if not before, no matter what input
$A\oplus B$ enumerating the cuts $L_\x$ and $R_\x$ is used.
%, $\Gamma$
%will enumerate $u_{i,s}$ into its upper cut and $v_{i,s}$ into its lower cut.
With $v<v_{i,s}$ and $u_{i,s}<u$, this makes it clear that $\Gamma$
enumerates our original (arbitrary) $u$ and $v$ into the correct cuts.
Therefore, on all inputs $A\oplus B$ enumerating $L_\x$ and $R_\x$,
$\Gamma$ does compute $\g(\x)=\f(\x)$ as required.
\qed\end{pf}

\section{Uniformity}
\label{sec:uniformity}

It was not the original purpose of this article to prove anything new.
The intention was to present a new proof of Lusin's Theorem
in real analysis, using known facts from computability theory
and descriptive set theory, and thus to illustrate and illuminate a connection
between the principles used in standard proofs of Lusin's Theorem
and the principles from computability which make our proof here work.
Nevertheless, certain uniformities and computability results became
apparent during the creation of the proof in Section \ref{sec:Lusin},
and in the end we have
an effective version of Lusin's Theorem.  Sometimes new ideas entail
new results, even when not intended to do so.

The substantial uniformity in the creation of $\g$ from $\f$ in our
proof of Theorem \ref{thm:Lusin} yields the function $h$ that we describe here.
Similar uniform versions hold for our simpler results in Theorems
\ref{thm:genericCantor}, \ref{thm:LebesgueCantor}, and \ref{thm:genericR}.

\begin{thm}
\label{thm:uniformLusin}
There is a computable total function $h:\Q\times\N\to\N$ such that,
for each fixed $\alpha$ and $S$ and each $S$-decidable presentation $\A$ of $\alpha$,
whenever an $\alpha$-jump-computable function $\f$ is given by the oracle computation
$$ \Phi_e^{(S\oplus A\oplus B)^{(\A)}}$$
for all enumerations $A\oplus B$ of the cuts of each $\x\in\R$, the function $\g(\x)$ defined by
$$ \Phi_{h(\ep,e)}^{E(\A)\oplus S^{(\A)}\oplus A\oplus B}$$
realizes Lusin's Theorem \ref{thm:Lusin} for this $\f$ and an arbitrary rational $\ep>0$.
Here $E(\A)$ is the elementary diagram of $\A$, given as a subset of $\N$ by a G\"odel coding.

That is, the Turing functional for computing $\g$ can be determined uniformly
from that for $\f$, uniformly in $\ep$ and independently of the level of $\f$
in the Borel hierarchy up to $\alpha$.  
%(Of course, the oracle $S^{(\A)}$ used to compute $\g$ does depend on the Borel level of $\f$, as well as on $S$ itself.)
\qed\end{thm}

By the Padding Lemma (see e.g.\ \cite[Lemma I.3.2]{S87}), $h$ may also be assumed injective.
Of course, $h$ only describes how to determine the program for computing this $\g$.
In order to compute $\g$, one also needs the oracle $S^{(\A)}$ and the elementary
diagram $E(\A)$ of the linear order $\A$.  However, the choice of program for $\g$
depends only on $\ep$ and the program given for computing $\f$, not on
which $S$ and $\A$ are used to compute $\f$.  In fact, not all of the atomic diagram
$E(\A)$ is required as an oracle; it suffices to know which elements of $(\A,\prec)$ are limit points
from the left, which is the zero element, which pairs $(m,n)$ are adjacencies
(with $m\prec n$ and no elements between them), and whether $\A$ itself
is a limit ordinal, a successor, or zero.  For uniformity, though,
even the finite information, such as knowing which element is the left end point of $\A$,
must be given.  Finally, it is not necessary to be given this diagram itself:
since an $S$-oracle is given, one only needs to know an index for computing
this information about $\A$ from $S$.

Nowhere in the proof of Theorem \ref{thm:Lusin} did we use the fact
that the values $\f(\x)$ were finite real numbers (as opposed to $\pm\infty$).
Indeed, the same proof would work even if enumerations of improper Dedekind
cuts were allowed as outputs.  Moreover, the construction in Theorem \ref{thm:Lusin}
ensures that the points $\x$ where $\g(\x)\neq\f(\x)$ all have $\g(\x)$ finite,
as they all lie in $\U$, which is the union of intervals on which $\g$ is defined
to be finite.

\begin{thm}
\label{thm:infLusin}
For every Borel-measurable function $\f:\R\to\R\cup\{\pm\infty\}$ and every $\ep>0$,
there exists a continuous function $\g:\R\to\R\cup\{\pm\infty\}$ such that $\mu(\set{\x\in\R}{\f(\x)\neq\g(\x)})\text{$~\!<\ep$}$.  Moreover, the value $\g(\x)$ is infinite only when
$\f(\x)$ is, and all uniformities described in Theorem \ref{thm:uniformLusin} still hold here.
\qed\end{thm}

\section{The Causes of Discontinuity}

Our proof of Theorem \ref{thm:Lusin} emphasizes a remarkable fact about Lusin's Theorem.
Once the parameters $\ep$, $S$, and $\A$ are fixed, the proof uses
the exact same error set $\U_{\ep,S,\A}$ for every
function $\f$ it is given.  So one may legitimately argue that the non-continuity of functions $\f$
at this level of the Borel hierarchy is the ``fault'' of the real numbers $\x$ in
$$\U_{S,\A} = \bigcap_{\ep>0}~ \U_{\ep,S,\A},$$
namely, the set of those $\x\in\R$ that are not generalized-$\alpha$-low relative to $S$
(using the presentation $\A$ of $\alpha$).
In this sense, the functions themselves are not the obstacle:  their non-continuity was
caused by our inability to approximate $(S\oplus L_{\x}\oplus R_{\x})^{(\A)}$
for those $\x$ in the error intervals.

\begin{prop}
\label{prop:Lusinacrossfcts}
Fix $\alpha$, $S$, an $S$-decidable presentation of $\alpha$, and a rational $\ep>0$.
Then for every $\alpha$-jump $S$-computable function $\f:\R\to\R$,
the procedure in Theorem \ref{thm:uniformLusin} produces a continuous $\g$
such that the set $\set{\x\in\R}{\f(\x)\neq\g(\x)}$ is always contained
within the same open set $\U_{\ep,S,\A}\subseteq\R$
of measure $<\ep$, independent of the choice of $\f$.  Indeed, for $\ep_0\leq\ep_1$,
we have $\U_{\ep_0,S,\A}\subseteq\U_{\ep_1,S,\A}$.
\qed\end{prop}

The immediate objection to this proposition is that, just by translating $\f$ by a certain fixed
parameter $\bf{c}$, one could define a function $\f_{\c}(\x)=\f(\x-\bf{c})$ for which most of the
discontinuities of $\f$ move out of $\U_{\ep,S,\A}$.  This is true, but it requires $\bf{c}$ to be noncomputable,
indeed not $S$-computable, and so $\f_{\c}$ does not belong to the class of functions
considered in Proposition \ref{prop:Lusinacrossfcts}.  In fact, $\U_{\ep,S,\A}$
is closed under translation by $S^{(\A)}$-computable parameters, and under other similar gambits
one might concoct.  

The more informed objection to the proposition is that it is obvious:  there are
only countably many $\alpha$-jump $S$-computable functions, so by applying
Lusin's Theorem to the $n$-th such function with tolerance $\frac{\ep}{2^{n+1}}$,
we immediately prove the proposition.  This is correct, but the spirit of the proposition
is that it was not necessary to slice up the $\ep$-amount of measure this way:
our proof of Theorem \ref{thm:Lusin} defined $\U_{\ep,S,\A}$ using basic
computability theory, and then uniformly constructed some continuous $\g$ for each $\f$
such that they differed only within $\U_{\ep,S,\A}$.  Probably the best way to express this
is to note that the restriction of every such $\f$ to the complement of each
$\U_{\ep,S,\A}$ is itself $S^{(\A)}$-computable and hence continuous on this domain,
and that each such domain is simply a $\Pi^{S^{(\A)}}_1$ set of real numbers.

In contrast, however, the restriction of such an $\f$ to the complement of
$\U_{S,\A}$ (defined just above) need not be continuous.  Analogously,
while only measure-$0$-many real numbers fail to be generalized-$\alpha$-low
relative to $S$, no single Turing functional can compute $(S\oplus X)^{(\alpha)}$
from $S^{(\alpha)}\oplus X$ for all but measure-$0$-many $X$.
As an example of a $2$-jump-computable function $\f$ such that no restriction of $\f$
to a set $\D\subseteq\R$ of full measure is continuous on the domain $\D$, consider
$$ \f(\x) = \sum_{e\in (L_\x\oplus R_\x)'} \frac1{2^e}.$$
The set $(L_\x\oplus R_\x)'$ is c.e.\ relative to $(L_\x\oplus R_x)$, hence
uniformly decidable using the oracle $(A\oplus B)''$ whenever $(A\oplus B)$
enumerates the cuts of $\x$.  Thus $\f$ is $2$-jump-computable.  However,
from any enumeration of the cuts of $\f(\x)$, one can compute $(L_\x\oplus R_\x)'$,
uniformly in $\x$ and in the choice of enumeration.  Therefore,
by the last part of Lemma \ref{lemma:GL1}, the restriction of $\f$
to an arbitrary set of full measure cannot possibly be an $S$-computable
function, no matter what oracle set $S$ one chooses.  Then one
adapts Theorem \ref{thm:Weihrauch} to functions on subdomains within $\R$
to show that this restriction cannot be continuous.

\section{Computing Continuous Functions}
\label{sec:continuous}

When Lusin's Theorem is applied to a function $\f$ that is already continuous,
it holds trivially:  just take $\g=\f$.  One might ask whether the procedure
given in Section \ref{sec:Lusin} reflects this.  The immediate answer is that it does not:
if $\f$ is continuous but is presented to us as an $\alpha$-jump-computable function,
applying the procedure there will often produce a $\g$ that, while satisfying the requirements
of Lusin's Theorem, is not in fact equal to $\f$, not even up to a set of measure $0$.
For future investigation, we conjecture that this is inherent:
no uniform procedure (as in Theorem \ref{thm:uniformLusin})
instantiating Lusin's Theorem can also succeed in making $\g=\f$ when $\f$ is itself
continuous.

However, if we ask the same question restricted entirely to continuous functions $\f:\R\to\R$,
then it is possible to produce a procedure for computing the function from a procedure
for $\alpha$-jump-computing it.  (In general a stronger oracle is required, though.)
This situation could plausibly arise:
for example, perhaps we can only determine a jump-computation
for a solution $\f$ to some differential equation under certain initial conditions,
although such an $\f$, being differentiable, must be continuous.

\begin{thm}
\label{thm:diffeq}
Let $\alpha$ be a countable ordinal and $\A$ an $S$-decidable presentation of $\alpha$.
Then there exists a computable total injective function $h:\N\to\N$
such that, whenever
$$ \f = \Phi_e^{\left((S\oplus A\oplus B)^{(\A)}\right)}:\R\to\R$$
is an $\alpha$-jump $S$-computation of a continuous $\f$, we have
a $0$-jump $S^{(\A+1)}$-computation
$$ \Phi_{h(e)}^{\left( S^{(\A+1)}\oplus A\oplus B\right)}=\f.$$
%is an $S^{(\A+1)}$-computation of the same function $\f$.
\end{thm}
Here $\A+1$ is the presentation of the ordinal $\alpha+1$ with
$\dom{\A+1}=\dom{\A}\cup\{ k\}$, where the number $k=\min (\N-\dom{\A})$
is adjoined to $\A$ as a new greatest element.
\begin{pf}
Where in Theorem \ref{thm:Lusin}, the rational numbers were a hindrance to
be handled by errer sets, here instead they serve as our guide.
For every $q\in\Q$, the left and right cuts $L_q$ and $R_q$ are
computable uniformly in $\Q$, so $(S\oplus L_q\oplus R_q)^{(\A)}$
is $S^{(\A)}$-computable, uniformly in $q$, and an
$(S\oplus L_q\oplus R_q)^{(\A+1)}$ oracle can decide the set
$$ D=\set{(a,b,u,v)\in\Q^4}{(\forall q\in [a,b])~v< \f\text{$(q)<u$}}.$$
The elements of $D$ are ``boxes'' $(a,b)\times (v,u)$ in $\R^2$
within which the graph of $\f$ (restricted to $(a,b)$) must lie.
Now for any $\x\in\R$ and any enumeration $A\oplus B$ of the cut of $\x$,
we get an $S^{(\A+1)}$-computable enumeration of
$$E_{\x}=\set{(u,v)\in\Q^2}{(\exists a\in p_1(A))(\exists b\in p_1(B))~(a,b,u,v)\in D}.$$
By continuity there are boxes in $D$ with $u-v$ arbitrarily small,
and so the projections $p_3$ and $p_4$ of $E_{\x}$ are the right
and left cuts $R_{\f(\x)}$ and $L_{\f(\x)}$.  Thus we have a computation
of $\f$ below an $S^{(\A+1)}$-oracle, whose program is uniform in the index $e$.
\qed\end{pf}

\section{Appendix: Computable functions on $\R$}
\label{sec:computableanalysis}

Turing computability normally applies to functions from $\N$, the set
of all nonnegative integers, into itself.  By fixing a computable bijection
between $\N$ and $\Q$, we may equally well consider functions $\Q\to\Q$,
or $\Q\to\N$, or $\N\to\Q$.  We also use a standard computable bijection
$\la m,n \ra$ mapping $\N\times\N$ onto $\N$, and similarly for $\Q$.

We write $a,b,q,r, u,v$ and sometimes $\ep$ for rational numbers, and $\a,\b,\r,\x,\y$ for real numbers.
Real numbers correspond bijectively to Dedekind cuts under their usual definition:
nonempty downward-closed proper subsets of $\Q$ with no greatest element.
For our purposes, this definition must be adapted slightly.
\begin{defn}
\label{defn:enumcut}
The \emph{Dedekind cut} of $\x\in\R$ is the pair $(L,R)$,
where $L=\set{q\in\Q}{q<\x}$ and $R=\set{q\in\Q}{q>\x}$.
In particular, if $\x\in\Q$, then $\x$ $\!\notin L\cup R$.
We also define \emph{generalized Dedekind cuts} to include the pairs
$(\emptyset,\Q)$ and $(\Q,\emptyset)$, corresponding to $-\infty$ and $+\infty$,
in addition to the \emph{proper} Dedekind cuts defined above.

An \emph{enumeration} of a generalized Dedekind cut $(L,R)$ is a set that,
when expressed as a join $A\oplus B$, satisfies $\pi_1(A)=L$ and $\pi_1(B)=R$,
where $\pi_1(\la q,n\ra)=q$ is the projection map.
The \emph{canonical enumeration} of $(L,R)$ is the join
$(L\oplus\N)\oplus(R\oplus\N)$.
\end{defn}

It is natural to regard $A$ and $B$ as subsets of $\Q\times\N$,
so that an oracle for $A\oplus B$ allows us to list out the elements in the
projections of $A$ and $B$, and thus to enumerate both the lower cut and
the upper cut of the real $\x$.

We can now give a definition
of computability for functions on $\R$ using Dedekind cuts, instead of the
usual fast-converging Cauchy sequences, to represent real numbers.

\begin{defn}
\label{defn:computablefunction}
For each subset $S\subseteq\N$, a function $\f:\R\to\R$ is \emph{$S$-computable}
if there exists a Turing functional $\Phi$ such that, whenever 
$X$ is an enumeration of the Dedekind cut of any $\x\in\R$,
$\Phi^{S\oplus X}$ is the characteristic function of an enumeration
of the Dedekind cut of $\f(\x)$.

It is equivalent, and often simplifies matters, to have $\Phi^{S\oplus X}$
be a partial function from $\Q$ into $\{ 0,1\}$, understanding
$\set{q}{\Phi^{S\oplus X}(q)\converges=0}$ and $\set{q}{\Phi^{S\oplus X}(q)\converges=1}$
to be the left and right Dedekind cuts of $\f(\x)$, respectively.  In this case,
if $\f(\x)$ itself is rational, then $\Phi^{S\oplus X}(\f(\x))$ never halts.
\end{defn}

The possibility that $L\cup R$ omits an element of $\Q$ is the reason
for considering enumerations of cuts.
If we had simply taken $L$ as an oracle, rather than an enumeration of $(L,R)$,
then for rational $q$, the characteristic
function of the interval $[q,+\infty)$ would have been computable;
similarly with $R$ and the interval $(-\infty,q]$.  On the other hand,
if we had required the cut of a rational $\y$ to include $\y$ itself
on one side or the other, then functions such as $\f(\x)=\x^2-2$
would not be computable.  (For that $\f$, given an enumeration of
the cuts of $\x$~$=\sqrt2$, $\Phi$ would never be able to place the
rational $0$ with certainty on either side of the cut of $\f(\sqrt2)$.)

Classically it is standard to express Definition \ref{defn:computablefunction}
using \emph{fast-converging Cauchy sequences}
instead of Dedekind cuts.  (A Cauchy sequence $\la q_n\ra_{n\in\N}$
\emph{converges fast} to $\x$ $\!=\lim_n q_n$ if, for every $n$, $|q_n-\x|<$ $2^{-n}$.)
Readers will understandably be baffled at first by our choice to use enumerations of Dedekind cuts
as inputs and outputs of functions, rather than following tradition.
We request forbearance:  in our view, this is the simplest way to prove
Theorems \ref{thm:genericR} and \ref{thm:Lusin}, because
(as suggested in Definition \ref{defn:enumcut}) it affords a canonical
enumeration of the Dedekind cuts of each $\x$, whereas in general
no Cauchy sequence converging fast to $\x$ is identified as the canonical
such sequence.  Were it not for the simplification of the proofs 
of these theorems, we would gladly use the traditional definition,
which has proven appropriate for all work so far in computable analysis.

Definition \ref{defn:enumcut} makes Definition 
\ref{defn:computablefunction} equivalent to the usual definition of computable functions
on $\R$.  The next lemma proves this, by showing that we can pass effectively
between the different methods of representing real numbers.
%It also covers another, very similar way of representing real numbers:
%Dedekind cuts $(L_\pi,R_\pi)$ in the set $\pi+\Q=\set{\pi+q}{q\in\Q}$.
%These are defined essentially as in Definition \ref{defn:enumcut}: the cut of $\x$ in $\pi+\Q$
%is $(\set{q\in\Q}{\pi+q<\x},$~$\!\set{q\in\Q}{\pi+q>\x})$, and will be used in Section \ref{sec:nearcontinuity}.

\begin{lemma}
\label{lemma:conversion}
There exist Turing functionals converting each of the following representations
of real numbers $\x$ into the other:
\begin{itemize}
\item
An arbitrary enumeration of a Dedekind cut for $\x$ in $\Q$.
%\item
%An arbitrary enumeration of a Dedekind cut for $\x$ in the set $\pi+\Q$.
\item
An arbitrary Cauchy sequence that converges fast to $\x$.
\end{itemize}
\end{lemma}
\begin{pf}
The lemma states that, for example, there exists a Turing functional $\Upsilon$
such that, whenever $A\oplus B$ is an enumeration of the Dedekind cuts in $\Q$
for some real number $\x$, $\Upsilon^{A\oplus B}$ will be a Cauchy sequence converging
fast to the same $\x$.  Indeed, given an enumeration $A\oplus B$ of the cut
of some $\x$, on input $n$, $\Upsilon$ searches for
$\la a,s\ra\in A$ and $\la b,t\ra\in B$ with $b-a<\frac1{2^{n-1}}$.  For the first such
pair of pairs to be found, it outputs $q_n=\frac{a+b}2$
as the $n$-th term.  Since $\x\in$ $\!(a,b)$, this $q_n$ is strictly within $\frac1{2^n}$ of $\x$,
so $\la q_n\ra_{n\in\N}$ converges fast to $\x$.  (This is the
definition of \emph{fast convergence}: $|q_n-\x|<\frac1{2^n}$.)
Conversely, given any Cauchy sequence $\la q_n\ra_{n\in\N}$ converging fast to $\x$, the set
$$\left\{\la a,m\ra\in\Q\times\N~:~a<q_m-\frac1{2^m}\right\}\oplus
\left\{\la b,n\ra\in\Q\times\N~:~b>q_n+\frac1{2^n}\right\}$$
enumerates the Dedekind cuts of $\x$.
\qed\end{pf}

With Lemma \ref{lemma:conversion} it is clear that the well-known theorem
of Weihrauch holds for the functions of Definition \ref{defn:computablefunction}.
\begin{thm}[Weihrauch; see \cite{W00}]
\label{thm:Weihrauch}
A function $\f:\R\to\R$ is continuous if and only if there exists a set $S\subseteq\N$
such that $\f$ is $S$-computable in the sense of Definition \ref{defn:computablefunction}.
\end{thm}

\section{Appendix: Approximating the Iterated Jump}
\label{sec:ordinaljumps}

In this appendix we describe known concepts and results that may
be unfamiliar to readers outside computability theory, and also
describe the precise version of the iterated jump used in this article.
By definition, the \emph{jump}, or \emph{Turing jump}, $A'$ of a set
$A\subseteq\N$ is the relativization of the Halting Problem to the set $A$:
$$ A' = \set{e\in\N}{\Phi_e^A(e)\text{~halts}}.$$
Here $\Phi_0,\Phi_1,\ldots$ is the standard enumeration of all Turing functionals,
that is, of all programs for Turing machines endowed with an arbitrary ``oracle set''
of natural numbers.  One writes $\Phi_e^A$ for the partial function from $\N$
into $\N$ computed by the $e$-th such program when using the set $A$
as its oracle.  We refer the reader to \cite[Chap.\ III]{S87} for full definitions.
The set $\emptyset'$ is essentially the Halting Problem itself, and just as $\emptyset'$
is not computable, so likewise $A$ is always strictly below $A'$ in Turing reducibility:
$A <_T  A'$, by which we mean $A\leq_T A'$ but $A'\not\leq_T A$.  (It is not completely obvious
that one can compute $A$ using an $A'$-oracle, but the proof is not difficult.)

The \emph{jump operator} is simply the function $A\mapsto A'$, sending
each subset of $\N$ to its jump, and thus mapping the power set
$\P(\N)$ into itself.  This map is injective but not surjective, and it
preserves Turing reductions but can fail to preserve non-reductions:
the implication
$$ A\leq_T B \implies A'\leq_T B'$$
always holds, but its converse can fail.
Since the jump operator maps $\P(\N)$ into itself, it is natural
to iterate it.  We write $A^{(k+1)}$ for the jump of $A^{(k)}$, with $A^{(0)}=A$,
thus defining all finite jumps of each set $A$.

However, this is not the end of the iteration.  For $\omega$, the first
infinite ordinal, the \emph{$\omega$-th jump}
$A^{(\omega)}$ of a set $A$ is a sort of union of all the preceding jumps of $A$:
$$  A^{(\omega)} = \set{\la n,k\ra\in\N^2}{n\in A^{(k)}}.$$
One views the $k$-th jump of $A$ (for each $k$) as being coded into
$A^{(\omega)}$ as the $k$-th column, under the usual computable bijection
from $\N^2$ onto $\N$ mapping the ordered pair $(n,k)\in\N^2$
to its code number $\la n,k\ra\in\N$.
Clearly no finite jump $A^{(k)}$ can compute $A^{(\omega)}$, since if it could,
it would then compute the $(k+1)$-st column $A^{(k+1)}$, contradicting
the fundamental property of the jump operator.  We view $A^{(\omega)}$
as a natural ``next'' jump, in the sense of ordinals, after all finite jumps have been built.
(To be clear: the set $\omega$ is actually just the set $\N$, but now viewed as an ordinal.)

Nor yet is this the end of the process.  One now continues through successor
ordinals as before, with $A^{(\omega+1)}=(A^{(\omega)})'$ and so on.
At subsequent limit ordinals $\lambda$, it is not always as obvious as for $\omega$
exactly how to define the $\lambda$-th jump, but we can do so if given
a computable presentation of the ordinal $\lambda$ -- that is, a computable
linear ordering $\prec$ of the domain $\N$ such that $(\N,\prec)$
is isomorphic to $\lambda$ as a linear order.  (Later we will assume that
in the order $\prec$, the successor and limit relations are computable
as well.)  Then we can define the $k$-th column of $A^{(\lambda)}$ to represent
the $\alpha$-th jump $A^{(\alpha)}$, where $k\in\N$ is mapped to $\alpha\in\lambda$
by the isomorphism from $(\N,\prec)$ onto $\lambda$.  This gives a reasonable
notion of $A^{(\lambda)}$, except that it depends on the choice of the presentation
$\prec$ of $\lambda$.  With proper use of ordinal notations, one can now define
the set $A^{(\alpha)}$ for all ordinals $\alpha$ with computable presentations, and,
although the actual set depends on the notation chosen. its Turing degree does not.

Church and Kleene knew that there must be a countable ordinal with no computable
presentation, and the least such ordinal is now known as $\oCK$.  Iterating the jump
to $A^{(\oCK)}$ and beyond requires presentations of noncomputable ordinals.
In the work in this article, we are generally able to use as oracle a (noncomputable)
set $S$ capable of giving presentations of noncomputable countable ordinals.
Of course, for each $S$ there is a least countable ordinal $\omega_1^S$ which
has no $S$-computable presentation (and it then follows that no ordinal $>\omega_1^S$
has any $S$-computable presentation either).  On the other hand, every countable
ordinal is isomorphic to some linear order on $\N$, just by the definition
of countability, and so, for every countable $\alpha$, there is some oracle set
$S$ such that $\alpha<\omega_1^S$, i.e., such that $\alpha$ has an
$S$-computable presentation.  (Again, we will need a slightly stronger presentation
of $\alpha$, with successors and limits known, but with the right oracle this
can also be assumed.  In fact, Spector showed that $S$ itself can give such a presentation.)

Thus we have (sketchily) described the iterated notion of the jump operator.
Having done so, we will now give formal definitions for use in this article.
\begin{defn}
\label{defn:strong pres}
A \emph{presentation} $\A$ of a nonzero ordinal $\alpha<\omega_1$ is
a linear order $\A=(D,\prec)$ isomorphic to $(\alpha,\in)$, whose
domain $D$ is a coinfinite subset of $\N$ and whose least element
is the number $0$ in $\N$.  The presentation $\A$ is
\emph{$S$-decidable} if $S$ can compute the complete diagram of $\A$.
\end{defn}
A presentation is normally called
\emph{$S$-computable} if $S$ can compute its atomic diagram.
We will need more than just that here, and we could be more precise (as above)
about the exact information we require $S$ to compute:  the successor function
on $\A$, the existence of a greatest element in
$\A$, and the unary relations on $\A$ of having no immediate predecessor
and of being the least or the greatest element of $\A$.  Demanding that $S$
compute the complete diagram is overkill, but keeps the definition simple.
The requirement that $D$ be coinfinite is unusual in computable
structure theory, but helpful for our purposes here:
we will need it when $\alpha$ is a successor, in order to have a location
in which to code one more jump.

It is important to notice that, for every nonzero $k\in D$, the substructure $\B_k$ with domain
$\set{j\in D}{j\prec k}$ is a presentation of an ordinal $\beta_k <\alpha$, and that each $\beta<\alpha$
is isomorphic to $\beta_k$ for some unique $k\in D$.  Moreover, $S$ can compute the
complete diagram of each $\B_k$, uniformly in $k$.

\begin{defn}
\label{defn:alphajump}
For a presentation $\A$ of a nonzero ordinal $\alpha<\omega_1$, the \emph{$\A$-jump}
$C^{(\A)}$ of a set $C\subseteq\N$ is the subset of $\N$ containing those codes
$\la k,n\ra$ for pairs $(k,n)$ such that either:
\begin{itemize}
\item
$k\in\dom{\A}$ and $n\in C^{(\B_k)}$; or
\item
$\A$ has a $\prec$-greatest element $j$ and $k$ is the least number $>j$
in $\overline{D}$ and
$$ (\forall n\in\N)~[\la k,n\ra\in C^{(\A)} \iff \Phi_e^{\left( C^{(\B_j)}\right)}(e)\converges].$$
\end{itemize}
The ordinal $\alpha=0$ has only one presentation $\A_0$, with empty domain,
and we define $C^{(\A_0)}=\N\times C$ for every set $C\subseteq\N$.
\end{defn}
So, for $k\in D$, the $k$-th column of $C^{(\A)}$ is simply the set $C^{(\B_k)}$,
meaning that for every $\beta<\alpha$, the $\beta$-th jump of $C$ (under
the presentation $\B_k$) appears as the $k$-th column.  If $\alpha$ is a limit ordinal,
all other columns are empty, but for a successor ordinal $\alpha=\beta+1$, with $j$
the greatest element of $\A$, the $k$-th column of $C^{(\A)}$ (for
to the least number $k>j$ in $\overline{D}$) presents the $\alpha$-th jump
of $C$ (under the presentation $\A$).  With $D$ coinfinite, such a column does exist
and is found using the diagram of $\A$.

Notice that the $0$-th column $\set{n}{\la 0,n\ra\in C^{(\A)}}$ of $C^{(\A)}$
is just the set $C$ itself, for every presentation $\A$ of any ordinal,
because $0$ is always the least element of the presentation.  This gives us a uniform
way to recover $C$ from $C^{(\A)}$, independent of the presentation.

When $\alpha=m$ is finite, $C^{(\A)}$ is not literally the same set as the jump $C^{(m)}$ discussed
above, but for our purposes in computability they are equivalent:  one column of $C^{(\A)}$
actually is the set $C^{(m)}$, all others are computable from $C^{(m)}$, and we have a
$1$-reduction from the $k$-th column to $C^{(m)}$ uniformly in $k$.  More generally, for a presentation
$\A$ of a successor $\alpha=\beta+1$, one can view $C^{(\A)}$ as the jump of $C^{(\B_j)}$,
where $j$ is the greatest element in $\A$:  this essentially says that the $\alpha$-th jump
is the jump of the $\beta$-th jump.  ($C^{(\B_j)}$ itself also appears inside $C^{(\A)}$,
just as the $\beta$-th jump is $1$-reducible to the $\alpha$-th jump.)  The uniformity of Definition
\ref{defn:alphajump} across finite and infinite ordinals will simplify our arguments below.

The key property used in this article is that, for a fixed oracle set $S$,
``most'' sets $A$ satisfy $(S\oplus A)'\leq_T S'\oplus A$, and that this reduction
is uniform for most of those sets $A$.  It is impossible for this to hold for all $A$
(and in particular for $A=S'$), but it does hold for all $1$-generic sets $A$, as seen
the proof of Lemma \ref{lemma:generic} above.  For Lebesgue measure on $2^\N$,
Theorem 2 of \cite{S72} proves that a Turing reduction $A'\leq_T \emptyset'\oplus A$
exists for measure-$1$-many sets $A$.  It is not uniform on any set of full measure.
(That is, no single functional $\Phi$ suffices, even up to a set of measure $0$.) 

Since we wish to address functions on $\R$ as well as on $2^\N$, we need a more specific theorem,
using only (strict) Dedekind cuts $L\oplus R$ of real numbers as our oracles.  Since the set of all such
cuts has measure $0$ under Lebesgue measure on Cantor space (viewed here
as $2^{(\Q\oplus\Q)}$), we must re-prove the result of \cite{S72} for our own measure,
namely Lebesgue measure on $\R$.  

The first lemma is a warm-up for the main theorem, demonstrating the basic technique.
Fix a computable enumeration $q_0,q_1,\ldots$ of all rationals in the interval $(0,1)$.
We write $L_{\x}$ for the strict left Dedekind cut of a real number $\x$,
and $R_{\x}$ for its strict right cut.
For each $a=q_j$ and $b=q_k$ in $\Q$ such that $0\leq a<b\leq 1$ and every
$q_i\in (a,b)$ has $i>\max\{j,k\}$, we define the binary strings $\lambda_{a,b}$
and $\rho_{a,b}$ each to have length $l$, where $l$ is least with $q_l\in (a,b)$, and set
$$ \lambda_{a,b}(i)=\left\{\begin{array}{cl} 1, &\text{if~}q_i\leq a\\
0, &\text{if~}q_i\geq b\end{array}\right.~~~~~~~~\text{and~~}
\rho_{a,b}(i)=1-\lambda_{a,b}(i)=\left\{\begin{array}{cl} 1, &\text{if~}q_i\geq b\\
0, &\text{if~}q_i\leq a.\end{array}\right.$$
Thus the real numbers $\x$ in the open interval $(a,b)$
are precisely those $\x$ such that $\lambda_{a,b}$ is an initial segment of $L_{\x}$
(viewed as an infinite binary string) and $\rho_{a,b}$ is an initial segment of $R_{\x}$.
(These strings are examined further at the end of the proof of Theorem \ref{thm:uniformity}.)
For intervals $(a,b)$ where some $q_i\in (a,b)$ has $i<j$ or $i<k$, it is possible
to divide the interval into finitely many subintervals with the property required above,
effectively, and to consider a string $\lambda\oplus\rho$ for each subinterval.

%For completeness we include as possible values $a=-\infty$ and $b=+\infty$, with $\lambda_{a,b}$ defined exactly the same way.  Thus every initial segment of every $L_{\x}$ is of the form $\lambda_{a,b}$.

\begin{lemma}
\label{lemma:GL1}
For every rational $\ep>0$, there exists a Turing functional $\Psi_{\ep}$ such that,
for every $S\subseteq\N$,
$$ \mu(\set{\x\text{$~\!\in (0,1)$}}{\text{$\Psi_{\ep}^{S'\oplus L_{\x}\oplus R_{\x}}=(S\oplus L_{\x}\oplus R_{\x})'$}})
>1-\ep.$$
Moreover, there is a computable function $h$ mapping each $\ep>0$
to an index $h(\ep)\in\N$ such that $\Phi_{h(\ep)}=\Psi_{\ep}$.
So this process is \emph{uniform in $\ep$}, although when $\ep=0$, no $\Psi$ suffices.
\end{lemma}
\begin{pf}
Given an $e\in\N$, $\Psi_{\ep}^{S'\oplus L_{\x}\oplus R_{\x}}$
searches for a rational number $r\in [0,1]$, a finite initial segment $\sigma\subset S$, and a
finite collection $(a_0,b_0),\ldots,(a_m,b_m)$ of disjoint open rational subintervals of $(0,1)$
and a stage $s$ such that:
\begin{itemize}
\item
$(\forall i\leq m)~\Phi_{e,s}^{\sigma\oplus\lambda_{a_i,b_i}\oplus\rho_{a_i,b_i}}(e)\converges$; and
\item
$\sum_{i\leq m} (b_i-a_i) > r$; and
\item
there do \emph{not} exist a number $t$, a $\tau\subset S$, and finitely many disjoint rational intervals $(c_0,d_0),
\ldots,(c_n,d_n)$ within $(0,1)$ such that $\sum_{i\leq n} (d_i-c_i)\geq r+\frac{\ep}{2^{e+1}}$ and
$(\forall i\leq n)~\Phi_{e,t}^{\tau\oplus\lambda_{a_i,b_i}\oplus\rho_{a_i,b_i}}(e)\converges$.
\end{itemize}
The $S'$-oracle allows $\Psi_{\ep}$ to recognize the truth or falseness of the final
statement for any specific $r$, while the first two statements are decidable.
For the $r$ that is found, we have $r<\mu(\set{\y\text{$~\!\in (0,1)$}}{\text{$\Phi_e$}^{S\oplus
\text{$L$}_{\y}\oplus R_{\y}}\text{$(e)\converges$}}) \leq r+\frac{\ep}{2^{e+1}}$.  (This also makes it clear
why such an $r$ must exist:  arbitrarily much of the measure of this set
can be covered by finitely many initial segments of oracles $L_{\x}\oplus R_{\x}$.)
Now $\Psi_{\ep}$ examines the $L_{\x}\oplus R_{\x}$-portion of its oracle.
If, for some $i\leq m$, $\lambda_{a_i,b_i}\oplus\rho_{a_i,b_i}\subseteq L_{\x}\oplus R_{\x}$,
then it outputs $1$,
since such an $\x$ will lie in one of the intervals $(a_i,b_i)$.  For all other $\x$,
it outputs $0$, meaning that it thinks that $\Phi_e^{S\oplus L_{\x}\oplus R_{\x}}\text{$(e)\diverges$}$.
This output $0$ may be wrong for certain values $\x$, but only for $\frac{\ep}{2^{e+1}}$-many.
Since this holds for every $e$, incorrect outputs can only occur for at most
$\ep$-much of the interval $(0,1)$.
Moreover, it is clear that this procedure is uniform in $\ep>0$.
\qed\end{pf}

Lemma \ref{lemma:GL1} clearly can be repeated for all intervals $(n,n+1)$
with precision $\frac{\ep}{4^{|n|+1}}$, uniformly in $n\in\Z$, to make it hold on all of $\R$.
Finally, we will need to know that this lemma holds not only
for the jump operator, but for all iterations and relativizations of it.  Here is the full result.

\begin{thm}
\label{thm:uniformity}
There exists a Turing functional $\Psi_{\ep}$, uniform in the rational $\ep>0$,
such that for every set $S\subseteq\N$ and every fixed $S$-decidable
presentation $\A$ of any ordinal $\alpha<\omega_1^S$
(with $\alpha$-jumps $C^{(\A)}$ defined using this presentation),
the ``error set''
$$ \U_{\ep,S, \A}=\set{\x\text{$~\!\in [0,1]$}}{\text{$\Psi_\ep^{S^{(\A)}\oplus L_{\x}\oplus R_{\x}}\neq
(S\oplus L_{\x}\oplus R_{\x})^{(\A)}$}}$$
has measure $<\ep$ and is an $S^{(\A)}$-effective union of rational open intervals,
uniformly in $\ep$, $S$, $\alpha$, and $\A$.
Moreover, for all $\x$, $\Psi_\ep^{S^{(\A)}\oplus L_{\x}\oplus R_{\x}}$ is total.
%Finally, if we fix $S$ and $\alpha$ and a strong $S$-presentation $\A$ of $\alpha$,
%the program for $\Psi_\ep$ and the $S^{(\A)}$-computable enumeration of
%$\U_{\ep,S, \A}$ may be found uniformly in $\ep$.
\end{thm}

\begin{pf}
We give the procedure of $\Psi_{\ep}$ on an input $\la k,e\ra$, using the presentation
$\A$ of $\alpha$, whose complete diagram $\Psi_{\ep}$ can compute using its $S$-oracle.

For elements $k\in D=\dom{\A}$, the program runs in a highly recursive manner,
computing the $k$-th column of its output using (finitely much information from)
those columns whose numbers $i$ satisfy $i\prec k$ (according to the diagram of $\A$).
Since $\A$ is well-ordered by $\prec$, this procedure is well-founded and will eventually halt.
On input $\la k,e\ra$ with $k\in D$, the program checks whether $k$ is $0$ or
a limit point in $\A$.  For $k=0$, it uses its oracle to decide whether
$e\in S\oplus L_{\x}\oplus R_{\x}$ and outputs the answer.  For a limit point $k$, it decodes
$e=\la k',e'\ra$ and runs itself on this pair, since $\la\la k',e'\ra,k\ra\in (S\oplus L_{\x}\oplus R_{\x})^{\text{$(\A)$}}$
just if $\la k',e'\ra\in (S\oplus L_{\x}\oplus R_{\x})^{\text{$(\A)$}}$.  If $k$ is a successor, then
the program finds the immediate predecessor $i$ of $k$ under $\prec$,
using the diagram of $\A$, and attempts to determine whether $\Phi_e^I(e)\converges$,
where $I$ is the $i$-th column of the program's own output.  This requires running
the program itself many times, recursively, but only on finitely many inputs
(and only on pairs $\la i',e'\ra$ with $i'\preceq i$).  The procedure is the same
as in Lemma \ref{lemma:GL1}:  the program runs until it has found a finite
set of some measure $r$ of initial segments $\lambda_{a,b}\oplus\rho_{a,b}$
that will cause $\Phi_e(e)$ with oracle $\Psi_{\ep}^{S^{(\B_i)}\oplus\lambda_{a,b}\oplus\rho_{a,b}}$ to halt,
and has been told by the oracle $S^{(\A)}$ that the oracles of this form
that cause $\Phi_e(e)$ to halt have total measure at most $r+\frac{\ep}{2^{k+e+2}}$.  
Here we are using $\Psi_{\ep}^{S^{(\B_i)}\oplus\lambda_{a,b}\oplus\rho_{a,b}}$
as an approximation to $(S\oplus L_{\x}\oplus R_{\x})\text{$^{(\B_i)}$}$,
which is the actual content of the $i$-th column $I$.
The approximation is not always correct; below we will consider the measure
of the set on which it is incorrect, but for now the important point is that
it does always give an output, instead of diverging.  Finally, $\Psi_{\ep}$
determines whether $L_{\x}\oplus R_{\x}$ begins with any of the finitely many strings
$\lambda_{a,b}\oplus\rho_{a,b}$ that were found to make
$\Phi_e^{\Psi_{\ep}^{S^{(\B_i)}\oplus\lambda_{a,b}\oplus\rho_{a,b}}}(e)$ halt.
If so, it outputs $1$, while if not, it outputs $0$.

For elements $k\notin D$, the program checks whether any $j<k$
is the greatest element of $\A$, and if so, whether $k$ is the $<$-least
number in $\overline{D}$ greater than that $j$.
If not, then it immediately outputs $0$.  If so, then it runs in a similar
manner to the program in Lemma \ref{lemma:GL1}, recursively using the column
$J=\set{\la j,n\ra}{n\in\N}$ of its own output.  For definiteness we specify that
on input $\la k,e\ra$ (with $k\notin D$) it uses a tolerance of $\frac{\ep}{2^{k+e+2}}$
to approximate $J'$.

Since $\prec$ well-orders $\A$, it is readily seen (by induction on columns whose numbers
lie in $D$, ordered according to $\prec$) that this program halts on every input $\la k,n\ra$
with $k\in D$.  If $\alpha$ is a successor, the same proof then applies to the $k\notin D$
determined above, and for all other $k\notin D$ it halted immediately.
So the program $\Psi_{\ep}$ with arbitrary oracle $S^{(\A)}\oplus L_{\x}\oplus R_{\x}$
always computes a total function.  Next we consider the set $\U_{\ep,S,\A}$ of those
$L_{\x}$ on which it fails to compute $(S\oplus L_{\x}\oplus R_{\x})^{(\A)}$.  This can happen in many
ways.  For the very first jump, when $k_1$ is the second-to-left point of $\A$, the computation
on input $\la k,e\ra$ will be incorrect on a set of $\x$ of measure $<\frac{\ep}{2^{2+k_1+e}}$,
and so the set of those $\x$ for which there is an error anywhere in this column
has measure $< \frac{\ep}{2^{k_1+1}}$.  For the second jump, in column number $k_2$,
there are two reasons the computation could be incorrect:  either $\x$ lies in the set of measure
$<\sum_e \frac{\ep}{2^{2+k_2+e}} = \frac{\ep}{2^{1+k_2}}$ on which the approximation
goes wrong, or else the approximation was using an incorrect version of
$(S\oplus L_{\x}\oplus R_{\x})^{\text{$(\B_{k_1})$}}$ from column $k_1$.  However, we already counted
those $\x$ for which the $k_1$ column was incorrect, so the reals $\x$ added to the set
$\U_{\ep,S,\A}$ on account of column $k_2$ have total measure $<\frac{\ep}{2^{1+k_2}}$.
Similarly, for every $\x$ in $\U_{\ep,S,\A}$, either there is some $\prec$-least $k\in D$
such that the computation for $\x$ goes wrong in column number $k$, or else
$\alpha$ is a successor and the computation went wrong in column $k=\min(\overline{D})$.
Therefore, the total measure of $\U_{\ep,S,\A}$ is at most
$$  \left(\sum_{k\in D} \frac{\ep}{2^{1+k}}\right)+\frac{\ep}{2^{1+\min(\overline{D})}}
~~\leq~~\sum_{k\in\N}\frac{\ep}{2^{1+k}}~~=~~\ep.$$

The foregoing paragraph already essentially explained how we can uniformly enumerate
the open set $\U_{\ep,S,\A}$ from an $S^{(\A)}$-oracle.  Those $\x$
in $\U_{\ep,S,\A}$ for which the first column $k_1$ was incorrect all have
$\Psi_{\ep}^{S^{(\A)}\oplus L_{\x}\oplus R_{\x}}(e)=0$ for some $e$ such that eventually
$\Phi_e^{(S\oplus L_{\x}\oplus R_{\x})}\text{$(e)$}$ halted.  With the $S^{(\A)}$-oracle
we can run both of these computations with arbitrary strings of the form
$\lambda_{a,b}\oplus\rho_{a,b}$ in place of $L_{\x}$.  When we find any $\lambda_{a,b}\oplus\rho_{a,b}$
and $e$ for which $\Psi_{\ep}(e)$ gave $0$ but $\Phi_e^{S\oplus\lambda_{a,b}\oplus\rho_{a,b}}(e)\converges$,
we enumerate the open interval $(a,b)$ of $\R$ into $\U_{\ep,S,\A}$
(since all $\x$ there have $\lambda_{a,b}\oplus\rho_{a,b}\subseteq L_{\x}\oplus R_{\x}$).
For column $k_2$, we do the same, using $\Psi_{\ep}^{S'\oplus L_{\x}}$ to compute
the oracle for the computation $\Phi_e^ { (S\oplus L_{\x}\oplus R_{\x})^{(\B_{k_1})}}(e)$;
if it does so incorrectly, then this $\x$ was already enumerated at the previous step,
while if it does so correctly, then we will enumerate $\x$ into $\U_{\ep,S,\A}$ just if
$\Phi_e^ { (S\oplus L_{\x}\oplus R_{\x})^{(\B_{k_1})}}(e)$ halts and
$\Psi_{\ep}^{(S^{(\A)}\oplus L_{\x}\oplus R_{\x})^{(\B_{k_1})}}(e)=0$.  Likewise, every
$\x\text{$~\!\in\U_{\ep,S,\A}$}$ will eventually be enumerated by this process,
because some finite initial segment $\lambda\oplus\rho\subseteq L_{\x}\oplus R_{\x}$ must have
been adequate to cause both of these events to occur, and we will eventually find that segment
(which must be of the form $\lambda_{a,b}\oplus\rho_{a,b}$, with $a=\max(\lambda^{-1}(1)\cup\{0\})$
and $b=\min(\lambda^{-1}(0)\cup\{1\})$) and enumerate $\x$ into $\U_{\ep,S,\A}$.

Finally we discuss the situation
of isolated points in $\U_{\ep,S,\A}$.  The strings $\lambda_{a,b}\oplus\rho_{a,b}$
(and their substrings) are not the only possible initial segments of oracles $L_{\x}\oplus R_{\x}$:
the other possiblity occurs when $\x$ itself is equal to the rational number $q_j$, in which case
$q_j\notin L_{\x}$ and $q_j\notin R_{\x}$.  Initial segments $\lambda\oplus\rho$
of such strings still satisfy the property $\min(\lambda^{-1}(0)) \leq \max (\rho^{-1}(0))$,
but they are allowed to have at most one $j<|\lambda|$ for which $\lambda(q_j)=\rho(q_j)=0$
(with $\rho(i)=1-\lambda(i)$ for all other $i$).  

Intuitively (and by definition), all rational numbers $q_j$ should lie in the error set $\U_{\ep,S,\A}$.
However, our enumeration so far could have omitted some rationals.
To include them (thus showing that $\U_{\ep,S,\A}$ is indeed open), notice that
for each $q_j$, we can effectively find the indices $c$ and $d$ of two other relevant functionals:
\begin{itemize}
\item
$\Phi_c^{S^{(\A)}\oplus L_{\x}\oplus R_{\x}}$ halts just if $L_{\x}$ contains
some rational $>q_j$ (that is, just if $\x>q_j$); and
\item
$\Phi_d^{S^{(\A)}\oplus L_{\x}\oplus R_{\x}}$ halts just if $R_{\x}$ contains
some rational $<q_j$ (that is, just if $\x<q_j$).
\end{itemize}
$\Phi_c$ will contribute an error interval to $\U_{\ep,S,\A}$ of the form $(q_j,q)$
for some rational $q>q_j$, and $\Phi_d$ will contribute one of the form $(r,q_j)$.
Therefore, including $q_j$ itself in $\U_{\ep,S,\A}$ keeps it open,
as now the entire interval $(r,q)$ is contained in $\U_{\ep,S,\A}$.  So, along with the
error intervals previously enumerated into $\U_{\ep,S,\A}$, we also enumerate
the interval $(r,q)$ defined here for this $q_j$, noting that $r$ and $q$
were computed effectively from $j$  We do the same for every other rational
$q_j$ in $\Q$ as well.  Thus $\U_{\ep,S,\A}$ is still $S^{(\A)}$-effectively open,
and the countably many new points do not change its measure.
\qed\end{pf}

\comment{
In Theorem \ref{thm:uniformity}, it would have been difficult to work with all of $\R$ at once,
but any interval $[a,b]$ of finite measure (with $S$-computable end points)
could have been used in place of $[0,1]$, and the argument would be uniform
in those end points.  However, now that the theorem is proven, it is an easy matter
to do it uniformly simultaneously for all $[n,n+1]$ with $n\in\Z$, and with
$\ep_n = \frac43\cdot\frac{\ep}{2^{2+|n|}}$, to do the same for all of $\R$.
}

As mentioned above, it is easy to repeat this process for all intervals $[n,n+1]$.
We record this as our final corollary.
\begin{cor}
\label{cor:uniformity}
There exists a Turing functional $\Psi_{\ep}$, uniform in the rational $\ep>0$,
such that for every set $S\subseteq\N$ and every fixed $S$-decidable
presentation $\A$ of any ordinal $\alpha<\omega_1^S$
(with $\alpha$-jumps $C^{(\A)}$ defined using this presentation),
the ``error set''
$$ \U_{\ep,S, \A}=\set{\x\text{$~\!\in\R$}}{\text{$\Psi_\ep^{S^{(\A)}\oplus L_{\x}\oplus R_{\x}}\neq
(S\oplus L_{\x}\oplus R_{\x})^{(\A)}$}}$$
has measure $<\ep$ and is an $S^{(\A)}$-effective union of rational open intervals,
uniformly in $\ep$, $S$, $\alpha$, and $\A$.
Moreover, for all $\x$, $\Psi_\ep^{S^{(\A)}\oplus L_{\x}\oplus R_{\x}}$ is total.
\qed\end{cor}

\end{document}